\documentclass[a4paper,11pt,english]{article}

\usepackage[utf8]{inputenc}
\usepackage[T1]{fontenc}
\usepackage{babel}
\usepackage{amsfonts,amssymb,amsmath,amsthm}
\usepackage{dsfont,upgreek,eufrak}

\usepackage{fullpage}
\usepackage{units}

\usepackage{graphicx}
\usepackage{multirow}
\usepackage[all]{xy}
\usepackage{float}

\usepackage{color}
\definecolor{ao}{rgb}{0.0, 0.5, 0.0}

\usepackage{hyperref}
\hypersetup{
	colorlinks=true,
	linkcolor=blue,
	bookmarks=true,
	bookmarksopen=true,
	citecolor=blue,
	allcolors=blue
}

\usepackage{cite}
\usepackage{url}

\makeatletter
\newcommand{\xRightarrow}[2][]{\ext@arrow 0359\Rightarrowfill@{#1}{#2}}
\makeatother

\newtheorem{thm}{Theorem}[section]
\newtheorem{defi}{Definition}[section]
\newtheorem{rem}{Remark}[section]
\newtheorem{prop}{Proposition}[section]
\newtheorem{lem}{Lemma}[section]

\usepackage{bigints}
\usepackage{verbatim}
\usepackage[blocks]{authblk}

\author[1]{Duboux Thibaut}
\author[1]{Offret Yoann}
\affil[1]{Institut de Mathématiques de Bourgogne (IMB) - UMR CNRS 5584\\
	Université de Bourgogne Franche-Comté, 21000 Dijon, France}
\title{\bf \Large Maximum Entropy Random Walks: the Infinite Setting and the Example of Spider Networks with their Scaling Limits}

\begin{document}

\date{\today}
\maketitle

\noindent\rule{\linewidth}{1pt}

\vspace{5pt}

\noindent
{\bf Abstract.} In this article, we establish solid foundations for the study of Maximal Entropy Random Walks (MERWs) on infinite graphs. We introduce a generalized definition that extends the original concept, along with rigorous tools for handling this generalization. Unlike conventional simple random walks, which maximize entropy locally, MERWs maximize entropy globally along their paths, marking a significant paradigm shift and presenting substantial computational challenges. Originally introduced by physicists and computer scientists in \cite{ref19}, MERWs have connections to concepts such as Parry measures and Doob $h$-transforms. Our approach addresses the challenges of existence, uniqueness, and approximation, illustrated through examples and counterexamples. Even in the infinite setting, MERWs continue to maximize the entropy rate, albeit in a less direct manner. Additionally, we conduct an in-depth analysis of weighted spider networks, including scaling limits, revealing various phenomena characteristic of the infinite framework, notably a phase transition. A unified proof of scaling limits based on submartingale problems is presented. Furthermore, we consider some extended models, where the spider lattice provides valuable insights, highlighting the complexity of studying these walks for general infinite weighted graphs.

\vspace{10pt}
\noindent\rule{\linewidth}{1pt}

\vspace{10pt}

\vspace{10pt}


\vspace{10pt}

\tableofcontents
\vspace{10pt}

\section{Introduction}

The most popular way to randomly explore a locally finite graph $G$ without any additional information is to assume that a walker at a given node jumps to any neighboring node chosen uniformly at random, and does so independently at each time step. This stochastic process is known as a Simple Random Walk (SRW), or, as referred to in \cite{ref19}, a Generic Random Walk (GRW). This choice, among all possible random walks, can be justified by its property of maximizing entropy production at each step. Following many classical references, such as \cite{woess_2000}, we shall refer to a random walk on a graph simply as a Markov chain.\\

\noindent
\textbf{{Entropy rate of a random walk.}} The concept of entropy, introduced by Ludwig Boltzmann, is fundamental in the fields of Statistical Physics and Thermodynamics. Similarly, the field of Information Theory, developed by Claude Shannon in the 1940s, also recognizes the importance of this quantity. We refer to their groundbreaking papers \cite{Boltz,Shannon}.  Here, all we need to know is that the entropy of a distribution $ \mu $ on a countable set $ V $ is defined by $H(\mu) = -\sum_{x \in V} \mu(x) \ln(\mu(x))$. When $X$ is a random variable on $ V $, $ H(X) $ represents the entropy of the distribution of $ X $. Besides, if $ \text{card}(V) = N $ is finite, the maximum value of $ H(\mu) $ is attained when $ \mu $ is the uniform probability measure on $V$, and it equals $ \ln(N) $. Concerning Markov chains, the quantity of significant interest is the entropy rate $h$ (see \cite{khin,ref16} for instance). When  $ (X_n)_{n \geq 0} $ is an irreducible and positive recurrent Markov chain on $V$, $ h $ is independent of the initial distribution and depends only on the invariant probability measure $ \pi $ and the transition kernel $ P $ : 
\begin{equation}\label{rateent}
h \equiv \lim_{n\to\infty}\frac{H(X_0,\cdots,X_n)}{n} = -\sum_{x,y \in V} \pi(x) P(x,y) \ln(P(x,y)).
\end{equation}
As an example, for the GRW on a finite graph with vertex set $V$, one has 
\begin{equation}
\pi(x)=\frac{d(x)}{\sum_{y\in V} d(y)}, \quad P(x,y)=\frac{A(x,y)}{d(x)}, \quad\mbox{and}\quad  h_{\text{GRW}} = \frac{\sum_{x \in V} d(x) \ln(d(x))}{\sum_{x \in V} d(x)},
\end{equation}
where $A$ denotes the adjacency matrix of the graph, and $d(x)=\sum_{y\in V}A(x,y)$ represents the out-degree of the vertex $x$. \\

\noindent
\textbf{{A brief history of MERWs and their applications.}} Maximum Entropy Random Walks represent a paradigm shift from a local to a global perspective. In essence, these are random walks that maximize entropy along their paths or, equivalently, the entropy rate (\ref{rateent}). This approach was recently introduced in \cite{ref17,ref18,ref19}. Among their findings, the authors emphasize the strong localization phenomenon of MERWs in slightly disordered environments. This property is particularly relevant in Quantum Mechanics, especially in the context of the Anderson localization phenomenon (we refer to \cite{konig2} for a mathematical survey). More broadly, MERWs appear to hold significant implications for statistical physics (we can allude to \cite{TERNOVSKY,Pathcount}). The concept of MERW is closely related to that of Parry measures for subshifts of finite type, as defined in \cite{Parry} and recently explored in \cite{Marcovici}. They are also referred to as Ruelle-Bowen random walks (see \cite{Ruelle,RuelleBowen} for instance). This idea is also subtly present in \cite{Hetherington,ArnoldL} and an alternate interpretation of these random walks based on large deviation theory is given in \cite{Touchette}. Furthermore, MERWs could be instrumental in studying and modeling complex networks, as suggested in \cite{Sinatra,Delvenne,Demetrius}. Lastly, the MERW concept has found applications in diverse scientific areas, such as community detection \cite{OchabBurd,Ochab2}, link prediction \cite{Li}, and even quasispecies evolution \cite{Smerlak}.\\

\noindent
\textbf{{The finite setting.}}
While significant progress has been made in the mathematical framework of MERWs, further inquiry is yet required. Current advancements mostly pertain to finite graphs, which present advantageous properties. Specifically, when dealing with an irreducible finite graph $G$, the Perron-Frobenius theorem guarantees the existence and uniqueness of a MERW. As illustrated in \cite{ref17,ref18,ref19}, its Markov kernel $P$ and its invariant probability measure $\pi$, for all vertices $x,y$, can be written as
	\begin{equation}\label{kernel}
	P(x,y)=A(x,y)\frac{\psi(y)}{\rho\, \psi(x)}\quad\text{and}\quad \pi(x)=\varphi(x)\psi(x).
	\end{equation}
	Here, $A$ still denotes the adjacency matrix of the graph, $\rho$ is its spectral radius, and $\psi$ and $\varphi$ are respectively the positive right  and left $\rho$-eigenvectors of $A$, suitably normalized so that $\pi$ defines a probability measure. Besides, it can be shown that the corresponding entropy rate is $h_{\rm MERW}=\ln(\rho)$. Intriguingly, all trajectories of length $n$ between vertices $x$ and $y$ have the same probability, given by $
	\rho^{-n}{\psi(y)}/{\psi(x)}$.    
	While the trajectory distribution is not uniform, it becomes uniform when conditioned on trajectory length and endpoints. This property suggests the rich combinatorial features inherent in MERWs. Equation (\ref{kernel}) evokes the well-known Doob $h$-transform, commonly encountered when conditioning stochastic processes to remain within a specific domain. For relevant references, we refer to \cite{doob,Doney,konig1,Denisov,Denisov2} and Remarks \ref{Martin} and \ref{staypos}. However, we emphasize that the MERW perspective is quite different: there is no underlying random walk, the domain is fixed, and it is this domain that determines the probability transitions (see also Remark \ref{confu}). Furthermore, to broaden the scope, one can substitute the adjacency matrix $A$ with a weighted variant (strictly positive across edges) and require the MERW to maximize 
	\begin{equation}\label{rate0}
	h(Q)=-\sum_{x,y\in G}\mu(x)Q(x,y)\ln \left(\frac{Q(x,y)}{A(x,y)}\right),
	\end{equation}
	over the positive-recurrent Markov kernels $Q$ on $G$, where $\mu$ denotes the invariant probability  distribution of $Q$. When the entries  $A(x,y)$ are non-negative integers, this formulation  can be interpreted as a MERW on a multi-edge graph. Additional constraints, like energy conditions, can be introduced as discussed in \cite{ref17,dixit}. The positive eigenfunction $\psi$ is prominent when assessing node influence in complex networks, forming the crux of the eigenvector centrality method \cite{miserable}. For physicists, the function $\psi$ in (\ref{kernel}) can be interpreted as a wave function, specifically the ground state of the following discrete Schrödinger equation
	\begin{equation}
	-\Delta\psi(x) +H(x)\psi(x)=-\rho\,\psi(x),
	\end{equation}
	where $\Delta$ is the graph Laplacian and $H$ is the potential defined by
	\begin{equation}
	\Delta f(x)=\sum_{y\in G}A(x,y)(f(y)-f(x))\quad\text{and}\quad H(x)=-\sum_{y\in G}A(x,y).
	\end{equation}
	For symmetric matrices $A$, we have $\pi(x)=\psi^2(x)$: the stationary probability distribution of the MERW is the square of the wave function. There are only a limited number of solvable models where the spectral radius and the associated wave function are explicitly known and determining these in general is a challenging task. For specific examples, such as Cayley trees with a finite number of generation or periodic ladder graphs, we refer to \cite{Ochab} and \cite{ref18} respectively. Obviously, for sufficiently small graphs, it is feasible to compute these values numerically and carry out computer simulations of the MERW.
	
\begin{rem}\label{confu}
This method of generating random walks may initially seem confusing. Typically, transition probabilities are chosen in an ad hoc manner. Here, however, the primary object is the weighted graph, and the transition kernel of the MERW is merely a consequence of its structure. Once the network is fixed, the random walk to be studied is fully determined. For example, it is not possible, as is often the case, to impose additional conditions such as a null-drift criterion, square-integrable jumps, and so on.	
\end{rem}

\begin{figure}[H]
	\begin{center}
		\includegraphics[scale=1]{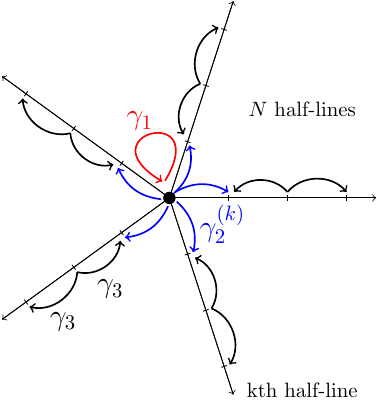}
		\caption{The weighted spider network}
		\label{model}
	\end{center}
\end{figure}

\noindent	
\textbf{Motivation, contribution and scope of this paper.} To our knowledge, there is a lack of consistent results for infinite networks. Some infinite periodic lattices are briefly investigated in \cite[Section 7.1]{ref18}, where some diffusion coefficients are computed, but no theoretical framework has been established. Additionally, some unweighted star graphs, as shown in Figure \ref{model}, are examined in \cite{TERNOVSKY,Pathcount}, but only from a combinatorial perspective. Phase transitions related to the number of paths are observed as the number of rays varies. However, no MERW is defined or investigated in these cases. Our main objectives are to go beyond the finite setting and begin expanding the bestiary of MERWs by providing compelling examples and counterexamples. In particular, we aim to address and elucidate the following questions:

\begin{itemize} 
\item[\textbf{{(a)}}] How can MERWs be properly defined on infinite graphs? Does a unique MERW exist? \item[\bf (b)] What about the entropy? What are the connections with MERWs on finite graphs? \item[\bf (c)] How do the scaling limits of MERWs compare with those of GRWs? 
\end{itemize} 
Additionally, we aim to clarify transversal issues: 
\begin{itemize} 
\item[\bf ($\ast$)] To which classical notions can this concept be related? What are the main challenges? What are the main tools to study these walks? 
\end{itemize}

We now detail the questions (a), (b), (c), and ($\ast$) by outlining the structure of this paper. First, we define the spectral radius $\rho$ of an infinite graph in Section \ref{sec:def}, based on the convergence parameter $R$ introduced by Vere-Jones \cite{ref30}. Computing $\rho$ can be challenging, as it requires asymptotic estimates of the number of walks. Next, analogous to the left-hand side of (\ref{kernel}), a MERW is defined using a positive eigenfunction $\psi$ associated with $\rho$. For infinite graphs, finding or approximating such a positive harmonic function can be a very difficult task.  Uniqueness and existence of MERWs are discussed in Section \ref{sec:exist}. In Section \ref{sec:reduc}, we discuss how graph symmetries can be used to reduce the computation of $\rho$ and $\psi$ to smaller, possibly finite lattices. In Section \ref{sec:MO}, we focus on the weighted spider network shown in Figure \ref{model}.
The investigation of this type of graph is motivated in \cite{TERNOVSKY,Pathcount} by the study of the conformational statistics of an ideal polymer chain. Here, the essentially unidirectional structure of this graph, coupled with nearest-neighbor interactions, allows for the explicit computation of the associated MERWs.
We demonstrate phase transitions between positive recurrent, null-recurrent, and transient behaviors, in accordance with those highlighted in \cite{TERNOVSKY,Pathcount} regarding the second moment asymptotics of the number of paths (the bifurcation region plays the role of an entropy trap). In this paper, the phase transitions find a more precise and quantified probabilistic interpretation, while introducing weights into the model. We also touch on non-nearest-neighbor and two-dimensional extensions in Section \ref{sec:extended}, illustrating how combinatorial problems can rapidly become more complex.

Definition \ref{MERW} may seem questionable for several reasons. First, it does not always ensure the maximization of the entropy rate (\ref{rate0}), as the corresponding Markov chain may not be positive recurrent. Additionally, if we replace $\rho$ with any $r > \rho$ and $\psi$ with an $r$-positive eigenfunction $\varphi$ in (\ref{def}), the resulting random walk still maximizes the entropy of paths of fixed length with given endpoints (see Theorem \ref{pruitt} and the example at the beginning of Section \ref{sec:entropy}). Why, then, choose $\rho$, and how does it relate to the entropy rate? Section \ref{sec:entropy} addresses this, synthesizing known results from dynamical systems in the context of MERWs. It turns out that $\ln(\rho)$ is the supremum of the entropy rate (\ref{rate0}) over all positive recurrent kernels $Q$, and even over all irreducible kernels on finite subgraphs (see Theorem \ref{entrate}). A surprising example, at the end of Section \ref{sec:entropy}, demonstrates that some MERWs, as defined in Definition \ref{MERW}, cannot be approximated by MERWs on finite subgraphs.\\

Regarding scaling limits, a classic example is Donsker's celebrated result \cite{Donsker}, which has inspired extensive research. Under certain conditions, Donsker's theorem shows that scaling limits of GRWs are Brownian motions. We aim to demonstrate that many continuous-time stochastic processes can be interpreted as scaling limits of MERWs. In Section \ref{sec:FSL}, we present the functional scaling limits derived  from the three type of MERWs obtained in Section \ref{sec:MO}. These limits include standard and drifted Walsh Brownian motions, as well as the three-dimensional Bessel
 process. The phase transitions presented above are also reflected in these scaling limits. In the more complex context of continuous-time processes, Section \ref{sec:continuous} shows how these diffusions minimize certain Kullback-Leibler divergences but in a less accessible manner. Section \ref{sec:proof} offers a unified proof of these limits using a submartingale problem approach, addressing challenges like the singularity of Walsh diffusions and the drift of the three-dimensional Bessel process at the origin. Notably, the MERW in the $R$-transient case on a spider lattice with $N=1$ is a Bessel-like random walk (or Lamperti Markov chain), whose scaling limit is studied in \cite{Lamp}. Finally, Section \ref{sec:extended} extends some of these results to particular two-dimensional networks and non-nearest-neighbor adjacency structures based on the spider lattice.

\section{General Framework}
\label{sec:GF}
\setcounter{equation}{0}

In what follows, let $G$ represent a countable irreducible weighted graph. We use $A$ to denote the weighted adjacency matrix and $\mathcal E$ for the set of edges. For simplicity, we will refer to $G$ as the set of vertices. Additionally, we shall assume that
\begin{equation}\label{bounddegree}
\sup_{x\in G} \sum_{y\in G}A(x,y)<\infty.
\end{equation}
When $G$ is unweighted, meaning that $A(x, y) \in \{0, 1\}$ for all $x, y \in G$, this condition simply means that the out-degrees of vertices are uniformly bounded.

\subsection{An Expanded Definition of MERWs}\label{sec:def}

For the primary results on infinite positive matrices that we make use of, we refer to \cite{ref30}.

\begin{defi}\label{radius}
	For any arbitrary $x,y\in G$, the combinatorial spectral radius, denoted by $\rho$,  is defined as the inverse of the radius of convergence for   
	$\sum_{n=0}^\infty A^{n}(x,y)z^n$. Notably, it is independent of the choice of $x$ and $y$.
\end{defi}

In essence, the leading asymptotic behavior of the number of $n$-step trajectories from $x$ to $y$ is on the order of $\rho^n$. Assumption (\ref{bounddegree}) above ensures that $\rho$ is finite.

\begin{defi}\label{MERW}
	A random walk on $G$ is termed a MERW if, for all vertices $x,y\in G$, its Markov kernel is defined as
	\begin{equation}\label{def}
	P(x,y)=A(x,y)\frac{\psi(y)}{\rho\,\psi(x)},
	\end{equation}
	where $\psi$ represents a positive eigenfunction of $A$ associated with the spectral radius $\rho$ (also referred to as a positive $\rho$-harmonic function).
	
\end{defi}

Analogously to (\ref{kernel}), if $\varphi$ is a positive left eigenfunction of $A$ associated with the eigenvalue $\rho$, then $\pi(x) = \varphi(x)\psi(x)$ is an invariant measure of the MERW. Besides, without loss of generality, we can assume that $\psi(o) = 1$ for a chosen base point $o \in G$. 

\begin{rem}\label{choquet}
The set of positive solutions $\varphi$ to $A\varphi = \rho \varphi$ with $h(\varphi) = 1$ is a convex compact set for the pointwise topology. Let $C$ denote this set and let $E\subset C$ be the subset of extremal solutions. By the Krein–Milman and Choquet theorems, for any $\psi\in C$, there exists a probability measure $\mu$ on $E$ such that  
\begin{equation}\label{martin0}
\psi(x) = \int_E \varphi(x)\, \mu(d\varphi).
\end{equation}

\end{rem}

\subsection{Existence and Uniqueness}
\label{sec:exist}

To ensure the existence and uniqueness of such MERW, we require further assumptions related to the recurrence and transience of Markov processes. The following definition is based on \cite{ref30} where $R=1/\rho$ is named the convergence parameter.

\begin{defi}\label{R-rec}
	Let $x, y \in G$ be arbitrary. The weighted adjacency matrix $A$ is termed ${R}$-recurrent (resp. $R$-transient) if 
	\begin{equation}
	\sum_{n\geq 0} \frac{A^n(x,y)}{\rho^n} = \infty \quad \left(\text{resp.} \sum_{n\geq 0} \frac{A^n(x,y)}{\rho^n} < \infty\right).
	\end{equation}
	If $A$ is $R$-recurrent, it is termed $R$-null (resp. $R$-positive) if $A^n(x,y)\rho^{-n}$ tends to zero (resp. does not tend to zero). Notably, these definitions are independent of the choice of $x$ and $y$.
\end{defi}

The following proposition follows easily from (\ref{def}), the preceding definition, and \cite{ref30}.

\begin{prop} \label{unique}
	Suppose that $A$ is $R$-recurrent. Then, $\rho$ is an eigenvalue of $A$, and there exist unique (up to a multiplicative constants) left and right eigenfunctions associated with $\rho$. Consequently, there is a unique MERW, which is recurrent. Moreover, this MERW is positive recurrent if and only if $A$ is $R$-positive. Furthermore, if $A$ is $R$-transient and a MERW exists, then it is necessarily transient.
\end{prop}

\begin{rem}\label{Martin}
In the $R$-transient situation, there may exist an infinite number of MERWs (see the example of spider lattice in Section \ref{sec:MO} for instance). Moreover, if $P$ is one of these Markov kernels, say associated with $\psi$ as in (\ref{MERW}), one can easily check that every positive harmonic function $h$ for $P$, that is, satisfying $Ph = h$, can be written as $h = \psi / \varphi$, where $\varphi$ is another positive solution of $A\varphi = \rho \varphi$. In particular, the classical Martin boundary theory can be used to describe all MERWs, similarly to (\ref{martin0}), and  each of these walks can be seen as a true Doob's $h$-transform of $P$ since 
\begin{equation}\label{doob}
P(x,y)\frac{h(y)}{h(x)}=A(x,y)\frac{\varphi(y)}{\rho\,\varphi(x)}.
\end{equation}
\end{rem}

\noindent
\textbf{An example of weighted graph with no existence.} In the $R$-transient case, neither the existence nor uniqueness of the MERW is assured. Consider  $G = \mathbb{N}_0 = \{0, 1, \cdots\}$ and  define $A$ such that $A(n, n-1) = 1$, $A(0, n) = \alpha_n>0$ for all $n \geq 1$, and $A(x, y) = 0$ elsewhere. It can be observed that
\begin{equation}\label{nolocally}
\sum_{n \geq 1} A^{n}(0, 0) z^n = \frac{1}{1 - \sum_{n \geq 1} \alpha_n z^{n+1}}.
\end{equation}
Therefore, if $\sum_{n \geq 1} \alpha_n z^{n+1}$ has a radius of convergence equal to $1$ and $\sum_{n \geq 1} \alpha_n < 1$,  we obtain that $\rho = 1$ and $A$ is $R$-transient. Since no solution exists for $\psi_0 = 1$, $\psi_n = \psi_{n-1}$ for all $n \geq 1$, and $\psi_0 = \sum_{n \geq 1} \alpha_n \psi_n$, no MERW exists for this setup.\\

As a matter of facts, it is possible to establish a necessary and sufficient condition for the existence of a MERW. This condition is closely tied to taboo-like probabilities and draws upon the foundational work by Harris and Veech on the existence of an invariant measure for a transient Markov chain. A detailed presentation of this result can be found in \cite{Pruitt}.

\begin{thm}\label{pruitt}
	The equation $A\psi=\lambda \psi$ with $\lambda>0$ has a positive solution $\psi$ if and only if one of the following conditions is satisfied: 
	\begin{enumerate}
		\item[(i)] $\lambda=\rho$ and $A$ is $R$-recurrent; 
		\item[(ii)] (a)  $\lambda=\rho$ and $A$ is $R$-transient, or (b) $\lambda>\rho$, and  in both cases, there exists an infinite subset $K\subset G$ and an exhaustive nested sequence $(G_j)_{j\geq 0}$ of $G$ with ${\rm card}(G_j)=j$ such that 
		\begin{equation}\label{HV}
		\lim_{j\to\infty,k\to\infty, k\in K}\frac{\sum_{y\notin G_j} A(x,y)\,{}_xF_{y,k}(\lambda^{-1})}{{}_xF_{x,k}(\lambda^{-1})}=0,
		\end{equation}
		where the power series ${}_xF_{y,k}(z)=\sum_{n=0}^\infty {}_xf_{y,k}^{(n)} z^{n}$ are recursively defined by
		\begin{equation}
		{}_x f_{y,k}^{(n+1)}=\sum_{w\neq x} A(y,w)\, {}_x f_{w,k}^{(n)}\quad\mbox{and}\quad {}_x f_{y,k}^{(0)}=\delta_{y,k}(1-\delta_{x,y}).
		\end{equation}
		Here $k\to\infty$ in $G$ in the sense of the Alexandroff extension, and  $\delta_{a,b}=1$ if $a=b$ and $\delta_{a,b}=0$ otherwise.
	\end{enumerate} 
\end{thm}

\begin{rem}\label{pruitt2}
	If $G$ is locally finite, then the Harris--Veech condition \eqref{HV} is satisfied. Consequently, there exists at least one MERW in this case. Note that the graph in the counterexample discussed around \eqref{nolocally} is not locally finite.
\end{rem}

\subsection{Entropy Rate characterizations}
\label{sec:entropy}

It should be noted that in case {\it (ii)(b)} of Theorem \ref{pruitt}, one can define a Markov kernel by replacing $\rho$ in (\ref{def}) with the corresponding $\lambda>\rho$. Conditionally on their length and their extremities, the probability of any trajectory remains proportional to its weight. One might question the reason for not replacing $\rho$ in Definition \ref{MERW} with an arbitrary $\lambda> \rho$ when feasible. The primary motivation is that we want MERWs to genuinely maximize the entropy production along the paths, in a manner yet to be defined.\\

\noindent
\textbf{A toy example.} Consider the standard lattice $G = \mathbb{Z}$. 
The set of extremal positive solutions to $\lambda \psi(x) = \psi(x+1) + \psi(x-1)$ with $\psi(0) = 1$ exists for any $\lambda = 2 \cosh(\alpha)$, $\alpha \geq 0$, and is given by $\psi_{\pm}(x) = e^{\pm \alpha x}$. 
The MERW $(X_n)_{n \geq 0}$ associated with $\psi_+$ is a usual biased random walk whose corresponding Markov kernel satisfies 
\begin{equation}
p_\alpha(x,x+1)=q_\alpha\quad\mbox{and}\quad p_\alpha(x,x-1)=1-q_\alpha,\quad\text{with } q_\alpha=\frac{e^{\alpha}}{e^{\alpha}+e^{-\alpha}}.
\end{equation}
For the one associated with $\psi_-$ it suffices to exchange $q_\alpha$ and $1-q_\alpha$. 
In any case, one has for all $x\in\mathbb Z$, 
\begin{align}
H(X_0,\cdots,X_n) & = \;    -\sum_{x_1,\cdots,x_n\in \mathbb Z} p_\alpha(x,x_1)\cdots p_\alpha(x_{n-1},x_n)\ln\left(p_\alpha(x,x_1)\cdots p_\alpha(x_{n-1},x_n)\right)\notag\\
& =\;  H(X_0,\cdots,X_{n-1}) - q_\alpha\ln\left(q_\alpha\right) - (1-q_\alpha)\ln\left(1-q_\alpha\right).\label{entropyrate}
\end{align}

It follows easily that, starting from an arbitrary point $X_0 = x$, the asymptotic rate of entropy $h$, defined by the first identity in (\ref{rateent}) (the second one being irrelevant in this case), is equal to the quantity involving $q_\alpha$ on the right-hand side of (\ref{entropyrate}), which is maximized for $\alpha = 0$, corresponding to $\lambda = \rho = 2$.\\

To delve deeper, recall that $h(Q)$ is defined in (\ref{rate0}) and introduce 
\begin{equation}\label{sup}
h^\star(G)=\sup\{h(Q) : \text{$Q$ is a positive-recurrent kernel on $G$}\}.
\end{equation}
Note that $h^\star(G)$ is bounded by the logarithm of the right-hand side of (\ref{bounddegree}). If $G$ is finite, the supremum of $h(Q)$ is attained at a unique positive recurrent kernel given by (\ref{kernel}), and we have $h^\star(G) = \ln(\rho(G))$, where $\rho(G)$ denotes the spectral radius of $G$. In the sequel, when considering a subgraph $H \subset G$, it is naturally endowed with the weight structure of $G$ through the restriction of the weighted adjacency matrix. The following result can be inferred from sources such as \cite{Seneta1,Seneta2,Seneta0,Ruette1,Ruette2,Sarig,GureSav,Sha,Shwartz1,Fayolle}. Further details will be provided below.

\begin{thm}\label{entrate}
	It holds that $h^\star(G)=\ln(\rho)$. Additionally, the supremum in (\ref{sup}) is actually a maximum if and only if $G$ is $R$-positive. When this condition is met, the maximum is attained by the unique MERW transition kernel. Moreover, one has
	\begin{equation}
	h^\star(G)=\sup\{h^\star(H) : \text{$H\subset G$ is finite and irreducible}\}.
	\end{equation}
	Furthermore, let $(H_n)$ be an exhaustive and increasing sequence of finite, irreducible subgraphs such that $h^\star(H_n)\longrightarrow h^\star(G)$ and let $P_n$ represent the unique MERW transition kernel on $H_n$. 
	\begin{enumerate}
		\item If $G$ is $R$-recurrent, then the sequence $(P_n)$ converges pointwise to the unique MERW transition kernel. 
		\item If $G$ is locally finite and $R$-transient, then the sequence $(P_n)$ is tight, and any of its limit points is a MERW transition kernel.   
	\end{enumerate}
\end{thm}

To be more specific, for unweighted graphs, $h^\star(G)$ represents the Gurevich entropy of the associated topological Markov chain. It has been established that an equilibrium measure (specifically, a Parry measure) exists if and only if $G$ is $R$-positive, and in such cases, this measure is unique. We refer the reader to \cite{Ruette1, Ruette2, Sha}. These results can be extended to weighted graphs using the concepts of topological pressure and potentials (see, for instance, \cite{Sarig, Shwartz1}). Regarding the convergence of a maximizing sequence, references \cite{Seneta1, Seneta2}, Chapter 6.4 of \cite{Seneta0} and \cite{Fayolle} provide insights, with the latter highlighting connections to the theory of large deviations. The reason for assuming that $G$ is locally finite in the context of transient graphs is due to the necessity of exchanging the limit and summation in the expression \begin{equation} \rho_n \psi_n(x) = \sum_{y \in H_n} A(x, y) \psi_n(y), \quad x, y \in H_n, \end{equation} where $\psi_n$ denotes the eigenfunction associated with the spectral radius of $H_n$, normalized such that $\psi_n(o) = 1$ for some fixed base point $o$ present in all $H_n$.

\begin{rem}\label{approx} At first glance, one might think that even in the transient case, all MERWs defined as in (\ref{MERW}) could be approximated by MERWs on finite subgraphs. However, as the following intriguing example demonstrates, this does not appear to be the case. 
\end{rem}

\noindent
\textbf{An example with quantized limit points.} Suppose $G = \mathbb{Z}$ carries the standard weight structure, except that $A(0, \pm 1) = \gamma$ with $\gamma > 0$. This is a specific case of the symmetric spider lattice investigated in Section \ref{sec:MO} and represented in Figure \ref{model} when $N=2$. In particular, we find that $G$ is $R$-transient if and only if $\gamma < 1$. Under these conditions, $\rho = 2$, and the two extremal eigenfunctions are
\begin{equation}
\psi^{(\pm)}(n)=
\left\{\begin{array}{ccl}
1+\Lambda n & \text{if} & \pm n\geq 0,\\
1 & \text{if} &  \pm n\leq 0,
\end{array}\right.\quad\text{with}\quad
\Lambda=\frac{2(1-\gamma)}{\gamma}.
\end{equation}
In particular, there is a one-to-one correspondence between $\{\lambda \psi^{(+)} + (1 - \lambda) \psi^{(-)} : 0 \leq \lambda \leq 1\}$ and the set of all MERWs. Let us introduce $H_{p,q} = \{-q+1, -q+2, \ldots, p-2, p-1\} \subset G$ for $p, q \geq 1$. Observe that $\rho(H_{p,q}) \uparrow 2$ as $p, q$ grow to infinity. Let $\psi_{p,q}$ denote the unique eigenfunction associated with the spectral radius $\rho(H_{p,q})$ and satisfying $\psi_{p,q}(0) = 1$.

\begin{prop} The set of all limit points of $\psi_{p,q}$ as $p, q \to \infty$ is given by the functions of the form $\mu_\delta \psi^{(+)} + (1 - \mu_\delta) \psi^{(-)}$, with $\delta \in \mathbb{Z} \sqcup {\pm \infty}$ and \begin{equation} \mu_\delta = \frac{\sqrt{\Lambda^2 \delta^2 + 4} + \Lambda \delta - 2}{2 \delta \Lambda} \in [0,1]. \end{equation} This probability is extended by continuity at $\delta \in {0, \pm \infty}$. In particular, the set of all MERWs that can be obtained as limits of classical MERWs on finite subgraphs is quantized. 
\end{prop}

\begin{proof}
Since $\rho(H_{p,q})<2$ we can write $\rho(H_{p,q})=2\cos(\theta_{p,q})$ where $\theta_{p,q}\in(0,\pi/2)$. It comes 
\begin{equation}
\psi_{p,q}(n)=\cos(\theta_{p,q}n)+(b_{-}\mathds 1_{n<0}+b_+\mathds 1_{n>0})\sin(\theta_{p,q}n),
\end{equation}
for some $b_\pm\in\mathbb R$. Analyzing the boundary conditions at points $0,p$, and $-q$, we derive $b_+-b_-=\Lambda{\rm cotan}(\theta_{p,q})$, 
$b_+=-{\rm cotan}(p\,\theta_{p,q})$ and
$b_-={\rm cotan}(q\,\theta_{p,q})$.
One can further express
\begin{equation}
\psi_{p,q}(n)=\frac{\sin((p-n)\,\theta_{p,q})}{\sin(p\,\theta_{p,q})}\mathds 1_{n \geq 0}+\frac{\sin((q+n)\,\theta_{p,q})}{\sin(q\,\theta_{p,q})}\mathds 1_{n < 0}.
\end{equation}
Since $\psi_{p,q}$ is positive, we obtain $\max(p,q)\theta_{p,q}<\pi$. Besides, since and $\theta_{p,q}\to 0$ as $p,q\to\infty$, one has $b^+-b^-\to \infty$ and then  ${\max(p,q)} \theta_{p,q}\to {\pi}$. 

Assume that $q \sim \alpha p$ for some $0 < \alpha < 1$. Then, we get  
$b_- \sim \cot(\alpha\pi)$
and 
\begin{equation}\label{asymp}
b_+ \sim \frac{1}{\pi - p\,\theta_{p,q}} \sim \frac{\Lambda p}{\pi}.
\end{equation}
From this, we infer that $\psi_{p,q} \to \psi^{+}$ pointwise, and subsequently, $\psi_{q,p} \to \psi^{-}$. As a result, we can identify the two extremal MERWs. 

Next, assume that $q \sim p$ with $p-q = \delta$ for some $\delta \in \mathbb{N}_0$. We get
\begin{equation}
\pi - q\,\theta_{p,q} \sim \left( 1 + \delta\frac{\theta_{p,q}}{\pi-p\,\theta_{p,q}}\right)(\pi-p\theta_{p,q}).
\end{equation}
We find that $b_{\pm} \sim \mu_{\pm} \Lambda$ where $\mu_- = 1-\mu_+$ and 
$\mu_+\left(1 + \frac{1}{1+\delta\mu_+\Lambda}\right) = 1.$ By symmetry, we deduce that the non-extremal MERWs obtainable through finite approximations are represented by the eigenfunctions $\mu_\delta \psi^{(+)} + (1-\mu_\delta)\psi^{-}$ for $\delta \in \mathbb{Z}$.
\end{proof}

\subsection{Automorphism and Reduced Models}
\label{sec:reduc}

In general, computing the combinatorial spectral radius and the associated eigenfunctions can be quite challenging. In this section, we provide tools to explore simpler models when the graph exhibits symmetries. For a deeper understanding of graph automorphisms and amenable groups, we refer the reader to \cite{woess_2000}. Let us define $\mathcal{T}$ as a subgroup of 
\begin{equation}
{\rm Aut}(G) = \left\{ \tau \in \mathfrak S(G) : \forall x,y\in G,\; A(\tau x,\tau y) = A(x,y)\right\},
\end{equation}
where  $\mathfrak S(G)$ denotes the symmetric group over $G$ and ${\rm Aut}(G)$ is the subgroup of graph automorphisms. The orbit of an element $x\in G$ under the action of $\mathcal{T}$ is represented by $\overline{x}$, and the entire orbit space is denoted by $\overline{G}$.

\begin{defi} \label{red}
	The set $\overline{G}$ is canonically endowed with a weighted graph structure inherited from that of $G$. The edges of this structure are defined as
	\begin{equation}
	\overline{\mathcal{E}} = \{(\overline{x}, \overline{y}) : \exists (x, y) \in \overline{x} \times \overline{y} \text{ such that } A(x, y) > 0\}.
	\end{equation}
	The corresponding weighted adjacency matrix is defined by $\overline{A}(\overline{x}, \overline{y}) = \sum_{y \in \overline{y}} A(x, y)$,	for any choice of $x \in \overline{x}$. Furthermore, if $G$ is irreducible, then $\overline{G}$ is also irreducible. This weighted graph is termed the {reduced graph}. In \cite{woess_2000}, it is also referred to as the factor graph $\mathcal{T} \setminus G$.
\end{defi}

If, for some $\lambda\in\mathbb{C}$ and $\overline{\psi} : \overline{G} \to \mathbb{C}$, we have $\overline{A}\,\overline{\psi} = \lambda \overline{\psi}$, then $A\psi = \lambda \psi$, where $\psi$ is defined by $\psi(x) := \overline{\psi}(\overline{x})$ for all $x \in G$. Conversely, if $A\psi = \lambda \psi$ and $\psi$ is $\mathcal{T}$-invariant, i.e., $\psi(\tau x) = \psi(x)$ for all $\tau \in \mathcal{T}$ and $x \in G$, then $\overline{A}\,\overline{\psi} = \lambda \overline{\psi}$, where $\overline{\psi}(\overline{x}) = \psi(x)$ for all $x \in G$.

\begin{rem}\label{lose}
	Unfortunately, eigenfunctions of $A$ are not necessarily $\mathcal{T}$-invariant, so finding all the $\lambda$-eigenfunctions of $\overline{A}$ does not guarantee that we have found all the $\lambda$-eigenfunctions of $A$.
\end{rem}

Furthermore, denote by $\overline{\rho}$ the combinatorial spectral radius of the reduced graph $\overline{G}$. Clearly, we have $\rho \leq \overline{\rho}$ because, for any $x, y \in G$ and $n \geq 1$, 
\begin{equation} \label{link}
\overline{A}^{\,n}(\overline{x}, \overline{y}) = \sum_{y \in \overline{y}} A^n(x, y).
\end{equation}
We shall provide conditions ensuring that $\overline{\rho} = \rho$. For a given $x_0 \in G$ and $n \geq 0$, let us define $B(x_0, n) = \left\{ x \in G : \exists\, 0 \leq k \leq n, A^k(x_0, x) > 0 \right\}.$ For any subset $L \subset G$, we define $\partial L$ as the set of vertices $y \in G \setminus L$ such that there exists $x \in L$ with $(x, y)$ being an edge of $G$. We use $|L|$ to represent the cardinality of $L$. Recall that $\mathcal{T}$ is termed quasi-transitive when $\overline{G}$ is finite.

\begin{prop}\label{reduction}
	We have $\rho = \overline{\rho}$ if any of the following conditions is met: 
	\begin{enumerate}
		\item[i)] There exists $x \in G$ such that $\overline{x}$ is finite.
		\item[ii)] $G$ is locally finite, $A$ is symmetric, and there exist $x_0,x \in G$ such that 
		\begin{equation}
		\lim_{n \to \infty} \frac{\ln(|B(x_0,n) \cap \overline{x}|)}{n} = 0.
		\end{equation}
		\item[iii)] There exists a positive $\mathcal{T}$-invariant function $\psi$ such that $A\psi \leq \rho \psi$.
		\item[iv)] $G$ is symmetric, locally finite, unweighted, $\mathcal{T}$ is quasi-transitive, and either
		\begin{enumerate}
			\item $G$ does not satisfy a strong isoperimetric inequality, i.e., 
			$\inf_{L \subset G, L \neq \emptyset} \frac{|\partial L|}{|L|} = 0;$
				\item or $\mathcal{T}$ is amenable and unimodular.
		\end{enumerate}
	\end{enumerate}
\end{prop}

\begin{proof}
	{\it $i)$} The power series $\sum_{n} A^n(x,y) z^n$, for $x,y \in G$, has a common radius of convergence $R = 1/\rho$ and possesses non-negative coefficients. From (\ref{link}), we conclude that $\overline{\rho} \leq \rho$.  
	
	{\it $ii)$} When $A$ is symmetric and locally finite, it can be viewed as a bounded linear operator on $\ell^2(G)$. Besides, it comes from \cite[Chap. II.10.]{woess_2000}  that   $\|A\|_2 = \rho$ and $\lim_{n \to \infty} \|A^n\|_2^{1/n} = \rho$. Specifically, considering $f_n(\cdot) = \mathds{1}_{B(x_0,n) \cap \overline{x}}$, we infer $\overline{\rho} \leq \rho$ from 
	\begin{equation}
	\limsup_{n \to \infty} \left( \overline{A}^{\, n}(\overline{x_0}, \overline{x}) \right)^{1/n} \leq \limsup_{n \to \infty} \|A^n f_n\|_2^{1/n} \leq \limsup_{n \to \infty} |B(x_0,n) \cap \overline{x}|^{1/2n} \|A^n\|_2^{1/n}.
	\end{equation}
	
	{\it $iii)$} By setting $\overline{\psi}(\overline{x}) = \psi(x)$, we get a positive function with the property $\overline{A}\,\overline{\psi} \leq \rho\,\overline{\psi}$ and thus we deduce that $\overline{\rho} \leq \rho$ by using  \cite{Pruitt}. 
	
	{\it $iv)$} Given that $\overline{G}$ is finite, there exists a positive function $\overline{\psi}$ such that $\overline{A}\,\overline{\psi} = \overline{\rho}\,\overline{\psi}$. Let $\psi$ represent the corresponding symmetric function on $G$ (a lift), and consider the random walk on $G$ with transition probabilities $P$ given by (\ref{def}), but with $\overline{\rho}$ in place of $\rho$. Given the symmetry of $A$, $\psi^2$ is a reversible measure which is bounded below and above. In other words, $P$ is a strongly reversible kernel. Referring to Theorem 10.3 and Corollary 12.12 in \cite{woess_2000}, we find that $\limsup_{n \to \infty} P^n(x,y)^{1/n} = 1$, which implies $\rho = \overline{\rho}$, whenever either condition $(a)$ or $(b)$ holds.
\end{proof}

\begin{rem}
	The preceding proposition can be applied to the infinite periodic lattices examined in \cite{ref18}, lending further rigor to their computation of the spectral radius.
\end{rem}

\section{Focus on spider MERWs}

\label{sec:MO}

\setcounter{equation}{0}
\subsection{Model and Settings}\label{reduced}

The model we consider (see Figure \ref{model}) is a star graph with $N$ half-lines perturbed at the origin.  It can be parameterized as ${G}=\big\{(n,k) :  n\in\mathbb N_0, k=1,\cdots,N\big\} \cup \{\mathbf 0\}$. For convenience, we make the identification $\mathbf 0=(0,1)=\cdots=(0,N)$. We will denote by $\partial G=\{\infty_1,\cdots,\infty_N\}$ the geometric boundary of $G$. Moreover, given any kernel $P(x,y)$ or function $\pi(x)$ on $G$, we denote by $P_k$ and $\pi_k$ their restrictions to the $k$-th leg $G_k$, and we often write $P_k(n,m)=P((n,k),(m,k))$ and $\pi_k(n)=\pi((n,k))$. The weighted adjacency matrix $A$  is defined for all $n\geq 1$ and $1\leq k\leq N$ by
$A(\mathbf 0,\mathbf 0)=\gamma_1$, $A_k(0,1)=\gamma_2^{(k)}$  and $A_k(n,n\pm 1)=\gamma_3$.  The tuple $\gamma = (\gamma_1, \gamma_2^{(1)}, \ldots, \gamma_2^{(N)}, \gamma_3)$ is assumed to belong to $\mathbb R_+\times (\mathbb R_+^*)^{N+1}$ and in the following, we shall set
\begin{equation}\label{lamb0}
\mathcal{S}_2 = \gamma_2^{(1)} + \cdots + \gamma_2^{(N)} \quad \text{and} \quad \Lambda = 2\gamma_3 - \gamma_1 - \mathcal{S}_2.
\end{equation}

\begin{defi}
	The case when $\Lambda=0$, $\Lambda<0$, or $\Lambda>0$ will be referred to as regular, attractive, or repulsive, respectively.
\end{defi}


\subsection{Spectral Radius}

We first observe that the spectral radius exhibits a phase transition phenomenon.

\begin{prop}\label{rho}
	The combinatorial spectral radius is given by
	\begin{equation}
	\rho=
	\begin{cases}
	2\gamma_3, & \text{if } \Lambda\geq 0,\\
	\dfrac{2\gamma_3(\gamma_1^2+\mathcal S_2^2)}{\gamma_1(2\gamma_3-\mathcal S_2)+\mathcal S_2\sqrt{\gamma_1^2+4\gamma_3(\mathcal S_2-\gamma_3)}}, & \text{if } \Lambda<0.
	\end{cases}
	\end{equation}
\end{prop}

\begin{proof} 
	Let $C_n$ be the $n$-th Catalan number. It is well-known that 
	\begin{equation} \label{G_S}
	S(z)=\sum_{n=0}^{\infty}C_n\gamma_3^{2n}z^{2n}=\frac{1-\sqrt{1-4z^2\gamma_3^2}}{2z^2\gamma_3^2}.
	\end{equation}
	The radius of convergence of $S(z)$ is $R_0= (2\gamma_3)^{-1}$. Let $R=\rho^{-1}$ be the radius of convergence of $F(z)=\sum_{n\geq 0}A^n({\bf 0},{\bf 0})z^n$. Using the classical arch-decomposition, see \cite[Chap. V.4.1]{flajoletSed} for instance,  one can write	
	\begin{equation}\label{G}
	F(z)=\frac{1}{1-(\gamma_1z +\mathcal S_2\gamma_3z^2S(z))}.
	\end{equation}

This standard method in algebraic combinatorics consists simply of decomposing excursions into elementary types. Here, we distinguish between those that start at $\mathbf{0}$ and return to $\mathbf{0}$ after one step and those that move from $\mathbf{0}$ to $1$ on some leg $k$, make an excursion from $1$ to $1$ while staying greater than $1$ on that leg, and then return to $\mathbf{0}$.

Thereafter, note that the function $x\mapsto \gamma_1x +\mathcal S_2\gamma_3 x^2 S(x) $ increases on $[0,R_0]$. Moreover, it can be verified that
	\begin{equation}\label{eq1}
	\gamma_1 R_0+\mathcal S_2\gamma_3 R_0^2 S(R_0)=\frac{\gamma_1+\mathcal S_2}{2\gamma_3}.
	\end{equation}
	It follows that $R=R_0$ when $\Lambda\geq 0$. If not, $R$ is the positive solution of  $\gamma_1 R+\mathcal S_2\gamma_3 R^2S(R)=1$, which is given by
	\begin{equation}\label{eq2}
	R=\frac{\gamma_1(2\gamma_3-\mathcal S_2)+\mathcal S_2\sqrt{\gamma_1^2+4\gamma_3(\mathcal S_2-\gamma_3)}}{2\gamma_3(\gamma_1^2+\mathcal S_2^2)}.
	\end{equation}
	This concludes the proof.
\end{proof}

\begin{rem}
	The spectral radius $\rho$ is identical to the model with a single leg where $\gamma_2^{(1)}=S_2$. When $\gamma_2^{(k)}\equiv \gamma_2$ is constant, this is a direct result of Proposition \ref{reduction}.
\end{rem}

\subsection{Markov Kernels}

\begin{prop}[regular case $\Lambda=0$]\label{reg_d}
	There exists a unique MERW. The positive right eigenfunction is given by $\psi\equiv 1$. For all $1\leq k\leq N$ and $n\geq 1$, the transition probabilities are
	\begin{equation}\label{P_reg}
	P_k(n,n+1)=\frac{1}{2},\quad P_k(n,n-1)=\frac{1}{2},\quad P_k(0,1)=\frac{\gamma_2^{(k)}}{2\gamma_3},\quad\text{and}\quad P(\mathbf 0,\mathbf 0)=\frac{\gamma_1}{2\gamma_3}.
	\end{equation}
	The process is null-recurrent with an invariant measure given by $\pi_k(n)=\gamma_2^{(k)}$ and $\pi(\mathbf 0)=\gamma_3$.
\end{prop}

\begin{proof}
	It can be readily verified that $\psi$ is a positive eigenfunction associated with $\rho=2\gamma_3$. Referring to the proof of Proposition \ref{rho}, it is apparent that $A$ is $R$-recurrent. Hence, Proposition \ref{unique} implies the uniqueness of the MERW. The remainder of the proof follows directly.
\end{proof}

\begin{prop}[attractive case $\Lambda<0$]\label{att_d}
	There exists a unique MERW. The positive right eigenfunction is $\psi_k(n)=\Gamma^n$ for all $n\geq 0$ and $1\leq k\leq N$. The factor $\Gamma$ is defined as
	\begin{equation}
	\Gamma = \frac{\rho-\gamma_1}{\mathcal S_2}.    
	\end{equation}
	The transition probabilities for all $1\leq k\leq n$ and $n\geq 1$ are given by
	\begin{equation}\label{P_att}
	P_k(n,n\pm 1) = \frac{\gamma_3}{\rho}\Gamma^{\pm 1},\ P_k(0,1) = \frac{\gamma_2^{(k)}\Gamma}{\rho},\ \text{and}\ P(\mathbf 0,\mathbf 0) = \frac{\gamma_1}{\rho}.
	\end{equation}
	Moreover, the MERW is positive recurrent. Its invariant probability measure is
	\begin{equation}\label{inv}
	\pi_k(n) = \frac{{\gamma_2^{(k)}}(1-\Gamma^2)\Gamma^{2n}}{{\mathcal S_2}\Gamma^2+\gamma_3(1-\Gamma^2)}\ \text{and}\ \pi(\mathbf 0) = \frac{{\gamma_3}(1-\Gamma^2)}{{\mathcal S_2}\Gamma^2+\gamma_3(1-\Gamma^2)}.
	\end{equation}
\end{prop}

\begin{proof}
	From Proposition \ref{rho}, we deduce the $R$-recurrence, leading us via Proposition \ref{unique} to the existence of a unique MERW on $G$. For all $n\geq 1$, consider  $\gamma_3\psi_k(n+1)+\gamma_3\psi_k(n-1) = \rho \psi_k(n)$ subject to
	\begin{equation} 
	\sum^N_{k=1}\gamma_2^{(k)}\psi_k(1) = \rho-\gamma_1\ \text{and}\ \psi(\mathbf 0) = 1.
	\end{equation}
	Let $\beta$ be the root of $\gamma_3 X^2 - \rho X + \gamma_3 = 0$ in the interval $(0,1)$. We can express $\psi_k(n)$ as
	$\psi_k(n) = a_k \beta^n + b_k \beta^{-n}$,
	with constants $a_k, b_k \in \mathbb{R}$. Using Proposition \ref{rho}, we find that
	\begin{equation}\label{relation}
	\gamma_1 + \mathcal S_2\Big[\gamma_3 \rho^{-1} S(\rho^{-1})\Big] = \rho,\quad\text{where}\quad 	S(z) = \frac{1-\sqrt{1-4z^2\gamma_3^2}}{2z^2\gamma_3^2}.
	\end{equation}
	From (\ref{relation}), we deduce
	\begin{equation}
	\beta = \frac{\rho-\sqrt{\rho^2-4\gamma_3^2}}{2\gamma_3} = \gamma_3 \rho^{-1} S(\rho^{-1}) = \Gamma.
	\end{equation}
	We can verify that the function $\psi_k(n) = \Gamma^n$ is indeed the unique solution. The invariant probability measure of the MERW is found by analyzing the left eigenvector of the system. For this, we solve the equation $\gamma_3\varphi_k(n+1) + \gamma_3\varphi_k(n-1) = \rho \varphi_k(n)$ for all $n\geq 2$ and $1\leq k\leq N$ subject to
	\begin{equation}
	\varphi_k(2) = \frac{\rho \varphi_k(1) - \gamma_2^{(k)}}{\gamma_3},\ \sum^N_{k=1}\varphi_k(1) = \frac{\rho-\gamma_1}{\gamma_3},\ \text{and}\ \varphi(\mathbf 0) = 1.
	\end{equation}
	We confirm that the function defined by $\varphi(\mathbf 0) = 1$ and $\varphi_k(n) = \frac{\gamma_2^{(k)}\Gamma^n}{\gamma_3}$ meets these conditions. Hence, it is the unique solution. The invariant probability measure is then obtained using standard computations.
\end{proof}

In the following, we set $\delta_{x,y}$ to be $1$ if $x=y$ and $0$ otherwise.

\begin{prop}[repulsive case $\Lambda>0$]\label{rep_d}
	There exists an infinite collection of MERWs generated by a finite number $N$ of linearly independent eigenfunctions $\{\psi^{(i)}: 1\leq i\leq N\}$. For all $1\leq k\leq N$ and $n\geq 0$, these are given by
	\begin{equation}\label{lamb}
	\psi_k^{(i)}(n) = 1 + \delta_{i,k}\frac{\Lambda}{\gamma_2^{(k)}}n.
	\end{equation}
	More precisely, there exists a one-to-one correspondence between MERWs and probability distributions $(\mu_i)_{1\leq i\leq N}$, through 
	\begin{equation}\label{psimu}
	\psi^{(\mu)} = \sum^N_{i=1} \mu_i \psi^{(i)}.
	\end{equation}
	The associated transition probabilities, for all $1\leq k\leq N$ and $n\geq 1$, are
	\begin{equation}\label{Q_rep}
	P_k^{(\mu)}(n,n\pm 1) = \frac{1}{2}\frac{\gamma_2^{(k)}+\mu_k\Lambda (n\pm 1)}{\gamma_2^{(k)}+\mu_k\Lambda n},\;
	P_k^{(\mu)}(0,1) = \frac{\gamma_2^{(k)}+\mu_k\Lambda}{2\gamma_3},\;
	P^{(\mu)}(\mathbf 0,\mathbf 0) = \frac{\gamma_1}{2\gamma_3}.
	\end{equation}
	Furthermore, let $\mathbb P^\mu_{x}$ denote the distribution of the MERW associated with $\mu$, starting from $x\in G$. Then, for all $1\leq k\leq N$, we have
	\begin{equation}\label{martinb}
	\mathbb P_0^\mu\left(\lim_{n\to\infty} X_n = \infty_k\right) = \mu_k.
	\end{equation}
\end{prop}

\begin{proof}
	We aim to solve $\psi_k(n+1) + \psi_k(n-1) = 2\psi_{k}(n)$ for all $1\leq k\leq N$ and $n\geq 1$, under the boundary conditions 
	\begin{equation}
\psi(\mathbf 0) = 1\quad\text{and}\quad \gamma_1 + \sum_{k=1}^N \gamma_2^{(k)}\psi_k(1) = 2\gamma_3.
	\end{equation}
It immediately follows that $\psi_k(n) = 1 + c_k n$ for some constants $c_k\geq 0$, which leads to the relation 
\begin{equation}
\sum_{k=1}^N \gamma_2^{(k)}c_k = \Lambda.
\end{equation}
This yields equations (\ref{lamb}), (\ref{psimu}), and (\ref{Q_rep}). Expanding upon this, we observe that $A$ is $R$-transient, meaning all the MERWs are transient. Let $(X_n)_{n\geq 0}$ be the MERW associated with the probability measure $\mu$. A harmonic function $h$ satisfies $\mathbb E_x[h(X_n)] = h(x)$ for all $x\in G$ and $n\geq 0$ if and only if $h = {\psi^{(\nu)}}/{\psi^{(\mu)}}$ for some other probability distribution $\nu$ (see Remark \ref{Martin}). The Martin boundary is thus represented by $\{1,\cdots,N\}$ and the Martin kernel is given by
	\begin{equation}
	K((n,k);i) = \frac{\psi^{(i)}_k(n)}{\psi^{(\mu)}_k(n)}.
	\end{equation} 
	Standard results on the Martin boundary of random walks assert that if $X_n$ starts from $\mathbf 0$, it almost surely converges within the Martin compactification to $i$ with probability $\mu_i$. The MERW corresponding to $\psi^{(i)}$ is  a classical symmetric nearest neighbor random walk in $\{(1,k), \cdots\}$ for all $k \neq i$. Given its transient nature, $\lim_{n\to\infty} X_n = \infty_i$ almost surely. Hence, we identify the Martin boundary with $\{\infty_1,\cdots,\infty_N\}$.
\end{proof}

\begin{rem}\label{staypos}
	Standard results (see \cite{Doney} for instance) indicate that the distribution of a simple symmetric random walk on $\mathbb{Z}$, conditioned to remain in $\mathbb{N}_0$, corresponds to the MERW on the spider lattice described in Proposition \ref{rep_d} when $N = 1$, $\gamma_1 = \gamma_2 = 0$, and $\gamma_3 = 1$ (a three-dimensional Bessel-like random walks, as investigated in \cite{Lamp,Kenn} and \cite[Chap. 3]{MPW}). However, MERWs cannot always be seen as non-trivial Doob's $h$-transforms of some other random walk, except when there is no uniqueness (see Remark \ref{Martin}).
\end{rem}

\section{Spider Functional Scaling Limits}
\label{sec:FSL}
\setcounter{equation}{0}

Let us introduce the space
\begin{equation}
\mathcal{G}=\big\{x=(\overline{x},k) : \overline{x}\in[0,\infty), k=1,\cdots,N\big\} \cup \{\mathbf{0}\}.
\end{equation}
Note that $G$ is canonically embedded in $\mathcal{G}$. We identify $\mathbf{0}=(0,1)=\cdots=(0,N)$ and denote by $\mathcal{G}_k=\{(\overline{x},k) : \overline{x}\geq 0\}$ the $k$-th leg. Furthermore, we equip $\mathcal{G}$ with the usual railway distance defined by 
\begin{equation}
d((x,i),(y,j))=|x-y|\delta_{i,j}+(x+y)(1-\delta_{i,j}).
\end{equation}
Restricted to $G$, this becomes the standard graph distance. For all $x=(\overline{x},k)\in\mathcal{G}$ and $\alpha\geq 0$, we set $\alpha x:=(\alpha\overline{x},k)$. For a proper planar embedding, the metric $d$ is equivalent to the usual Euclidean metric, and $\alpha x$ corresponds to the conventional scalar multiplication. Let $({\bf C},\mathcal{U})$ denote the space of continuous functions from $[0, \infty)$ to $\mathcal{G}$, equipped with the topology of uniform convergence on compact sets. We use $\Longrightarrow$ to signify the convergence in distribution of stochastic processes in $({\bf C},\mathcal{U})$ with the associated Borel $\sigma$-field. Let $\mathcal{F}_t$, $t\geq 0$, represent the canonical filtration on $\mathbf{C}$. For any sequence of real numbers $(X_n)_{n\geq 0}$, we define for all $t\geq 0$,
\begin{equation}
X_{t}=X_{\lfloor t\rfloor} + (t-\lfloor t\rfloor)(X_{\lfloor t\rfloor+1}-X_{\lfloor t\rfloor}).
\end{equation}
Here, $\lfloor x\rfloor$ denotes the largest integer less than or equal to $x$.

\begin{rem}
	It is possible to extend all the functional convergences discussed below to the space of {\it càdlàg} functions, either endowed with the usual Skorokhod topology or the uniform topology as described above. For more details, we refer to \cite[Chap. 18]{Billing}.
\end{rem}

\subsection{Regular Case}

We direct the reader to \cite{WalshYor} for the definition of the Walsh Brownian motion and to \cite{YorExcusrion} for the excursion theory of Brownian motion. Let $\{{\bf W}_{t}^{(\mu,x)}=(W_t,K_t) : t\geq 0\}$ be the standard Walsh Brownian motion on $\mathcal{G}$ starting from $x=(\overline{x},k)$ with spinning measure $\mu=(\mu_1,\cdots,\mu_N)$. Notably, when $W_t=0$, the value of $K_t\in\{1,\cdots,N\}$ is inconsequential. This process can be roughly described as follows. It is a continuous stochastic process on $\mathcal{G}$ where $W$ is a standard one-dimensional reflected Brownian motion starting from $\overline{x}$. It is noteworthy that $W_t$ can be expressed as $W_t=|B_t| = \mathcal{B}_t+L_t$, where $B$ and $\mathcal{B}$ are two standard one-dimensional Brownian motions starting from $\overline{x}$, and $L$ denotes the local time at 0 of $W$. To elaborate further, let $\tau$ represent the right-continuous inverse of $L$. Each excursion interval of $W$ away from zero can be expressed as $I_0=[0,\tau_{0})$ or $I_s=(\tau_{s-},\tau_s)$ for some $s>0$. The set difference of the union of these intervals is $\{t\geq 0 : W_t=0\}$ and has Lebesgue measure zero. Moreover, $K$ remains constant, say $\alpha_{I_s}$, over each $I_s$. We have $\alpha_{I_0}=k$ and, conditionally to $W$, $\{\alpha_{I_s} : s>0, I_s\neq \emptyset\}$ constitutes an independent collection of $\mu$-distributed random variables.

\begin{thm}[regular case $\Lambda=0$]\label{reg_c}
	Let $\{X_n\}_{n\geq 0}$ be the MERW presented in Proposition \ref{reg_d} and define 
	\begin{equation}\label{spin}
	\mu=\left(\frac{\gamma_2^{(1)}}{\mathcal{S}_2},\cdots,\frac{\gamma_2^{(N)}}{\mathcal{S}_2}\right).
	\end{equation}
	If $X_0=x_0$ is deterministic and depends on $L>0$ in such a way that for some $x\in\mathcal{G}$,
	\begin{equation}
	\frac{x_0}{\sqrt{L}}\xrightarrow[L\to\infty]{}x,
	\end{equation}
	then the following functional scaling limit holds:
	\begin{equation}\label{scale1}
	\left\{\frac{X_{Lt}}{\sqrt{L}}\right\}_{t\geq 0}\xRightarrow[L\to\infty]{}\{\mathbf{W}_t^{(\mu,x)}\}_{t\geq 0}.
	\end{equation}
\end{thm}

\subsection{Attractive Case}

The construction of the Walsh Brownian motion has been extended to various contexts. For Walsh diffusions, we direct the reader to \cite{FredWentz} for a functional analysis approach on graphs and to \cite{Walsh18,Walsh2,Walsh3,Walsh4,ref23} for semimartingale characterizations on rays. We allude to \cite{McKean} for the general Itô's theory of excursions. Fix $\lambda>0$ and $\overline x\geq 0$ and let $Z$ be the solution of the reflecting stochastic differential equation
\begin{equation}\label{reflect}
dZ_t = dB_t - \lambda dt + dL_t, \quad Z_t \geq 0, \quad Z_0 = \overline{x},
\end{equation}
with $B$ being a standard one-dimensional Brownian motion and $L$ a $B$-adapted, non-decreasing, continuous stochastic process which satisfies
\begin{equation}
\int_0^\infty \mathds{1}_{\{Z_t>0\}} dL_t = 0 \quad \text{and} \quad 
\int_0^\infty \mathds{1}_{\{Z_s=0\}} ds = 0 \quad \text{a.s.}
\end{equation}
Introduce the Walsh diffusion $\{\mathbf{Z}_t^{(\mu,x )}=(Z_t,K_t): t\geq 0\}$ on $\mathcal{G}$ starting from $x=(\overline{x},k)$ with the spinning measure $\mu$.  Similar to the Walsh Brownian motion, when $Z_t=0$, the specific value of $K_t\in\{1,\cdots,N\}$ is irrelevant. Moreover, we have $Z_t = |\mathcal{Z}_t|$, where $\mathcal{Z}$ is a (weak) solution of $d\mathcal{Z}_t = d\mathcal{B}_t - \lambda\, \text{sgn}(\mathcal{Z}_t) dt$, with $\mathcal{Z}_0 = \overline{x}$ and $\mathcal{B}$ a standard one-dimensional Brownian motion. Here $L$ represents the local time at zero for $Z$. The spinning measure $\mu$ is subject to the condition:
\begin{equation}
\forall i\in\{1,\cdots,N\}, \quad \lim_{\varepsilon\to 0^+}\frac{1}{2\varepsilon}\int_0^t \mathds{1}_{\{0< Z_s < \varepsilon\}} \mathds{1}_{\{K_s = i\}} ds = \mu_i L_t \quad \text{a.s.}
\end{equation}
As before, $K$ remains constant across each excursion interval $(I_s)_{s\geq 0}$ of $Z$. We have $K_t=k$ for $I_0$, and the values for $(I_s)_{s>0}$ are independent and distributed according to $\mu$, conditionally to $Z$. It is worth noting that $\mathcal{Z}$ is an ergodic diffusion with its reversible probability measure and one can check that the invariant probability measure of $Z$ is the exponential distribution of parameter $2\lambda$.

\begin{thm}[attractive case $\Lambda<0$]\label{att_c}
	Let $\{X_n\}_{n\geq 0}$ be the MERW described in Proposition \ref{att_d} with a $N+2$-tuple of parameters $\gamma=(\gamma_1,\gamma_2^{(1)}\cdots,\gamma_2^{(N)},\gamma_3)$ depending on $L>0$. Assume there exists $\zeta \in \mathbb R \times \mathbb R^{N} \times \mathbb R$ such that $\mathfrak Z_2 = \sum_{i=1}^N \zeta_2^{(i)} > 0$ and a positive constant $\lambda$ satisfying
	\begin{equation}\label{proj}
	\frac{\Lambda}{\mathcal S_2} \;\underset{L\to\infty}{\sim}\; -\frac{\lambda}{\sqrt L} \quad \text{and} \quad  \gamma - \zeta = \mathcal O\left(\frac{1}{\sqrt L}\right).
	\end{equation}
	Further assume that $X_0 = x_0$ is deterministic and that
	\begin{equation}
	\frac{x_0}{\sqrt L} \;\xrightarrow[L\to\infty]{}\; x.
	\end{equation}
	Then, for the spinning measure defined as $\displaystyle \mu = \left(\frac{\zeta_2^{(1)}}{\mathfrak Z_2}, \cdots, \frac{\zeta_2^{(N)}}{\mathfrak Z_2}\right)$, the following functional scaling limit holds:
	\begin{equation}\label{scale2}
	\left\{\frac{X_{Lt}}{\sqrt{L}}\right\}_{t\geq 0} \;\xRightarrow[L\to\infty]{} \; \{\mathbf Z_t^{(\mu,x) }\}_{t\geq 0}.
	\end{equation}
\end{thm}

\begin{rem}\label{heavytraffic} The stochastic process obtained when $N=1$ is the reflected Brownian motion with negative drift. This process plays a crucial role in queueing theory, particularly in the context of heavy-traffic approximations; see the seminal papers \cite{KingHeavy,KingHeavy2} and the survey \cite{heavywhitt}. In this context, the MERW can be interpreted as the workload of a discrete-time queueing system. Assumption (\ref{proj}) is essential to obtain a non-trivial limit, since the MERW converges to a stationary distribution by Proposition \ref{att_d}. \end{rem}

\subsection{Repulsive Case}

Firstly, introduce the well-known three-dimensional Bessel process $\{Y_t^{y}\}_{t\geq 0}$ starting from $y\geq 0$. This is the euclidian norm of a non-negative solution to the stochastic differential equation
\begin{equation}\label{bess3}
dY_t^{y}=dB_t+\frac{1}{Y_t^{y}}dt,\quad Y_0=y,
\end{equation}
where $B$ denotes a standard Brownian motion. For further details, we refer to \cite{Yor}. In essence, this is the euclidian norm of a $3$-dimensional Brownian motion. This is a transient Markov process satisfying $Y_t^{y}>0$ for all $t>0$, even when it starts at zero. 

Subsequently, for any $\overline x\geq 0$ and $1\leq k\leq N$, we consider the stochastic process $\big \{{\mathbf Y}_t^{(\overline x,k)}\big\}_{t\geq 0}$ on $\mathcal G$ defined by  
\begin{equation}
\mathbb P\left(\forall t\geq 0,\,{\mathbf Y}_t^{(\overline x,k)}=(Y_t^{\overline x},k)\right)=1.
\end{equation}
 This corresponds to the three-dimensional Bessel process on the $k$th leg. Furthermore, let $\mu=(\mu_k)_{1\leq k\leq N}$ be a probability distribution. We define the process $\big\{\mathbf Y_t^{(\mu,x)}\big\}_{t\geq 0}$ as
\begin{enumerate}
	\item For $x=\mathbf 0$: 
	$\mathbb P\left(\forall t\geq 0,\, {\mathbf Y}^{(\mu,\mathbf 0)}_t={\mathbf Y}_t^{(0,k)}\right)=\mu_k$ for all $1\leq k\leq N$.
	\item For $x=(\overline x,k)\neq \mathbf 0$:
	\begin{enumerate}
		\item If $\mu_k\neq 0$: $\mathbb P\left(\forall t\geq 0,\, \mathbf Y_t^{(\mu,x)}=\mathbf Y_t^{(\overline x,k)}\right)=1$.
		\item If $\mu_k= 0$:
		\begin{equation}
		\mathbf Y_t^{(\mu,x)}=\left\{\begin{array}{ll}
		(\overline x+W_t,k), & \text{for all } 0\leq t \leq \tau_0,\\
		{\mathbf Y}^{(\mu,\mathbf 0)}_{t-\tau_0}, & \text{for all } t\geq\tau_0,
		\end{array}\right.
		\end{equation}
		where $\tau_0=\inf \{t\geq 0\,|\,\overline x+W_t=0\}$ and $W$ is a standard one-dimensional Brownian motion, independent of $\mathbf Y^{(\mu,\mathbf 0)}$.
	\end{enumerate}
\end{enumerate}

\begin{thm}[repulsive case $\Lambda>0$]\label{rep_c}
	Let $\{X_n\}_{n\geq 0}$ be the MERW as specified in Proposition \ref{rep_d}, associated with the probability distribution $\mu$. If $X_0=x_0$ is deterministic and relates to $L>0$ such that
	\begin{equation}
	\frac{x_0}{\sqrt L}\;\xrightarrow[L\to\infty]{}\;x ,
	\end{equation}
	then the following functional scaling limits holds:
	\begin{equation}\label{scale3}
	\left\{\frac{X_{ Lt }}{\sqrt{L}}\right\}_{t\geq 0}\;\xRightarrow[L\to\infty]{} \; \{\mathbf Y_t^{(\mu,x)}\}_{t\geq 0}.	
	\end{equation}
\end{thm}

\begin{rem} Similarly to Remark \ref{staypos}, we note that a Brownian motion conditioned to remain positive is a three-dimensional Bessel process (see, for instance, \cite{Pitman}). Functional scaling limits of random walks conditioned to remain positive have been studied in \cite{Ingle0,Ingle1}. Furthermore, in \cite{Kenn}, a coupling between Bessel processes and Bessel-like random walks is constructed, providing an alternative proof of the scaling limit in this context. \end{rem}

\begin{rem} The aforementioned results can be extended  to the exclusion process involving two particles on $\mathbb{Z}$ that can jump left or right but cannot share a site. See \cite{Derrida,Schutz} for detailed reviews. By symmetry, $\psi(x,y) = 1 + (y - x)$ is a positive eigenfunction with spectral radius $\rho = 4$. In the scaling limit, one can obtained the equation 
	\begin{equation} d(Y_t - X_t) = dW_t + \frac{dt}{Y_t - X_t}, 
	\end{equation}
	where $X_t < Y_t$ are particle positions. Maximizing entropy reveals a standard electrostatic force. 
\end{rem}

\subsection{Continuous-Time counterparts of MERWs}
\label{sec:continuous}

In light of the scaling limits described, one may wonder whether the limit processes can be interpreted as maximal entropy stochastic processes without directly involving MERWs. 
In what follows we present a rather informal discussion, meant to illustrate the type of techniques that can be employed in the continuous setting. 
We focus on the case $N=1$ and explore the possibility of interpreting the three-dimensional Bessel process and the solution of (\ref{reflect}) as Maximal Entropy Stochastic Processes.\\

\noindent
\textbf{Kullback–Leibler Divergence (KLD).} Let $\gamma_n$ denote the (uniform) distribution of the first $n$-steps of the simple random walk on the regular graph $\mathbb Z$. It is noteworthy that maximizing the entropy $H(X_0,\cdots,X_n)$ -- n being fixed --  is equivalent to minimizing the Kullback–Leibler divergence (relative entropy) $D_{\rm KL}(\nu_n|| \gamma_n)$ where $\nu_n$ is the distribution of $(X_0,\cdots,X_n)$. 

Given two probability measures $\nu,\gamma$ where $\nu$ is absolutely continuous with respect to $\gamma$, the KL-divergence is defined as
\begin{equation}\label{KLD}
D_{\rm KL}(\nu\|\gamma)=\int \ln\left(\frac{d\nu}{d\gamma}\right) d\nu.
\end{equation}

To adapt this definition to continuous stochastic processes, we replace $\gamma$ with a reflected Brownian motion $W$, which satisfies 
$dW_t = dB_t + dL_t$, 
where $(B_t)_{t\geq 0}$ is a Brownian motion adapted to $(\mathcal F_t)_{t\geq 0}$ and $L_t$ denotes the local time of $W$ at $0$.
 The class of stochastic processes absolutely continuous with respect to $W$ will be defined as follows.

Let $\psi$ be a nonnegative, absolutely continuous function on $[0,\infty)$. 
Set $U=\{\psi>0\}$ and $\tau=\inf\{s\geq 0 : \psi(W_s)=0\}$ and define, for all $t\geq 0$, 
\begin{equation}\label{expmart}
M_t=\exp\left(\int_0^t\frac{\psi^\prime(W_s)}{\psi(W_s)}\,dB_s
-\frac{1}{2}\int_0^t \left(\frac{\psi^\prime(W_s)}{\psi(W_s)}\right)^2 ds\right)\mathds 1_{\{t<\tau\}}.
\end{equation}
For any $n\geq 1$, introduce the stopping times
\begin{equation}
\sigma_n := \inf\left\{t\geq 0 : \int_0^t \left(\frac{\psi^\prime(W_s)}{\psi(W_s)}\right)^2 ds \geq n\right\}\wedge n
\quad\text{and}\quad
\tau_n := \inf\Big\{t\geq 0 : \psi(W_t)\leq \tfrac1n\Big\},
\end{equation}
and set $\zeta_n := \sigma_n \wedge \tau_n$ and $\zeta := \sup_{n\geq 1}\zeta_n$. 
By Novikov’s condition (see for instance \cite{Yor}), the stopped process $(M_{t\wedge \zeta_n})_{t\geq 0}$ is a true martingale under $\mathbb P_x$ for every $x\in U$. 
Thereafter, one can define a probability measure $\mathbb Q_x^{(n)}$ on $\mathcal F=\sigma\big(\bigcup_{t\geq 0} \mathcal F_t\big)$ by setting, for all $t\geq 0$, 
\begin{equation}
\frac{d\mathbb Q_x^{(n)}}{d\mathbb P_x}\Big|_{\mathcal F_t} = M_{t\wedge \zeta_n}.
\end{equation}
Since $\mathbb Q_x^{(n+1)}$ and $\mathbb Q_x^{(n)}$ coincide on $\mathcal F_{\zeta_n}$, it follows that there exists a unique probability measure $\mathbb Q_x$ on $\mathcal F_\zeta := \sigma\big(\cup_{n\geq 1}\mathcal F_{\zeta_n}\big)$ such that 
$\mathbb Q_x|_{\mathcal F_{\zeta_n}} = \mathbb Q_x^{(n)}$ for all $n\geq 1$.

 Besides, applying Girsanov’s theorem (see again \cite{Yor}), under $\mathbb Q_x$ the process
\begin{equation}\label{newbrownien}
\widetilde B_t = B_t - \int_0^{t\wedge \zeta} \frac{\psi^\prime(W_s)}{\psi(W_s)}\,ds,
\end{equation}
is a Brownian motion up to the lifetime $\zeta$.
As a consequence, under $\mathbb Q_x$ the process $W$ satisfies the reflected stochastic differential equation
\begin{equation}
d W_t = d\widetilde B_t + \frac{\psi^\prime(W_t)}{\psi(W_t)}\,dt + d L_t, 
\quad 0\leq t<\zeta,
\end{equation}
with initial condition $ W_0 = x$, where $L$ still denotes the local time of $W$ at $0$.

Furthermore, assuming in addition
\begin{equation}\label{borne}
\int_U (\psi^\prime(x))^2\,dx<\infty,
\end{equation}
it follows from \cite[Theorem~6.3.4, p.~345]{Fuku} that $(W_t)_{t\geq 0}$ is, under $\mathbb Q_x$, a conservative diffusion, symmetric with respect to $\psi^2(x)\,dx$, which never hits the zero set of $\psi$ in finite time. In other words, $\tau=\zeta=\infty$ $\mathbb Q_x$-a.s.\@ for all $x\in U$. 
\begin{rem}
Although in \cite{Fuku} this result is proved for $\psi^2(x)\,dx$-almost every $x$ in a more general setting, in the present one-dimensional continuous case it extends to all $x\in U$ easily.	Besides, we recall that $(X_t)_{t \geq 0}$ is a $\psi^2(x) \, dx$-symmetric Markov process on $U$ if $\psi^2(x) \, dx$ is a reversible measure: for all sufficiently smooth test functions $f$, $g$ and for all $t \geq 0$,
\begin{equation}
\int_U P_t f(x) g(x) \psi^2(x) \, dx = \int_U f(x) P_t g(x) \psi^2(x) \, dx,
\end{equation}
where $(P_t)_{t \geq 0}$ is the Markov semigroup associated with $(X_t)_{t\geq 0}$. 
\end{rem}

Consequently, using \eqref{KLD} and \eqref{newbrownien} we obtain
\begin{align}\label{entropycont}
D_{\rm KL}\!\left(\mathbb Q_x|_{\mathcal F_{t\wedge \sigma_n}} \,\big\|\, \mathbb P_x|_{\mathcal F_{t\wedge \sigma_n}}\right)
&= \mathbb E_{\mathbb Q_x}\!\left[\int_0^{t\wedge \sigma_n} \frac{\psi^\prime(W_s)}{\psi(W_s)}\, dB_s 
- \frac{1}{2}\int_0^{t\wedge \sigma_n} \left(\frac{\psi^\prime(W_s)}{\psi(W_s)}\right)^2 ds \right] \\
&= \mathbb E_{\mathbb Q_x}\!\left[\int_0^{t\wedge \sigma_n} \frac{\psi^\prime(W_s)}{\psi(W_s)}\, d\widetilde B_s 
+ \frac{1}{2}\int_0^{t\wedge \sigma_n} \left(\frac{\psi^\prime(W_s)}{\psi(W_s)}\right)^2 ds \right] \label{eq:martinbtilde}\\
&=\frac{1}{2}\,\mathbb E_{\mathbb Q_x}\!\left[\int_0^{t\wedge \sigma_n} \left(\frac{\psi^\prime(W_s)}{\psi(W_s)}\right)^2 ds \right].\label{eq:div}
\end{align}
We need to localize up to $\sigma_n$ to ensure that the stochastic integral with respect to $\widetilde B$ in \eqref{eq:martinbtilde} is a true martingale under $\mathbb Q_x$ (hence has mean zero). Letting $n\to\infty$, since $\sigma_n\uparrow\infty$ as $n\to\infty$, the right-hand side of \eqref{eq:div} converges (by monotone convergence) to 
\begin{equation}
H_t=\frac{1}{2}\,\mathbb E_{\mathbb Q_x}\!\left[\int_0^{t} \left(\frac{\psi^\prime(W_s)}{\psi(W_s)}\right)^2 ds \right].
\end{equation}

Assuming that $\pi(dx)=\psi^2(x)\,dx$ is a probability measure, we obtain from the ergodic theorem that for any $x\in U$ and $s>0$, the asymptotic relative entropy rate  satisfies
\begin{equation}\label{casrecpos}
\lim_{t\to\infty}\frac{H_t}{t}=\frac{1}{2s}\,\mathbb E_{\mathbb Q_\pi}\!\left[\int_0^{t} \left(\frac{\psi^\prime(W_s)}{\psi(W_s)}\right)^2 ds \right]=
\frac{1}{2}\int_U \big(\psi^\prime(x)\big)^2\,dx.
\end{equation}

~

\noindent
\textbf{Repulsive case.} Assume that $\psi(x) > 0$ on $U=(0,L)$ and $\psi(x) = 0$ otherwise, for some $L > 0$. 
We are looking for a function that minimizes \eqref{casrecpos}. 
We obtain that, for every sufficiently small $\varepsilon \in \mathbb{R}$ and every sufficiently smooth function $h$ with compact support in $]0,L[$,
\begin{equation}\label{entropyh}
\int_{0}^L  \big(\psi^\prime(x) + \varepsilon h^\prime(x)\big)^2 \, dx 
+ \ell \left( \int_0^L \big(\psi(x) + \varepsilon h(x)\big)^2 \, dx - 1 \right) 
\geq \int_{0}^L \big(\psi^\prime(x)\big)^2 \, dx,
\end{equation}
where $\ell$ denotes a Lagrange multiplier. 
By examining the first-order term in $\varepsilon$, and integrating by parts, we obtain the equation $-\psi^{\prime\prime}(x) + \ell \psi(x) = 0$. 
By positivity of $\psi(x)$ and the boundary conditions, we easily find that the solution (whose square is a probability density) is given by
\begin{equation}
\ell = -\left(\frac{\pi}{L}\right)^2 
\quad \text{and} \quad 
\psi(x) = \sqrt{\tfrac{2}{L}} \, \sin\!\left(\tfrac{\pi}{L}x\right).
\end{equation}
Letting $L\to\infty$, we retrieve the drift of three-dimensional Bessel process,  since for all $x>0$, one has ${\psi^\prime(x)}/{\psi(x)} \sim {1}/{x}$.

~

\noindent
\textbf{Attractive case.} When $\gamma_1=1+\frac{\lambda}{\sqrt L}$, $\gamma_2^{(1)}=0$ and $\gamma_3=1$, we need to add another constraints on $\psi$ to retrieve the reflected Brownian motion with negative drift. We require that
\begin{equation}\label{mean}
\int_0^\infty x\psi^2(x)dx=\frac{1}{2\lambda}.
\end{equation}

Let $\psi$ be a local minimizer of \eqref{casrecpos}, and let $h$ be an arbitrary compactly supported smooth function with $h(0)=0$. As before, by considering $\psi+\varepsilon h$ for sufficiently small $\varepsilon\in\mathbb{R}$ and examining the first-order terms in $\varepsilon$, we obtain that necessarily  
\begin{equation}\label{EDO}
-\psi^{\prime\prime}(x)+\ell \psi(x)+ \kappa x\psi(x)=0,
\end{equation}
where $\ell$ and $\kappa$ are Lagrange multipliers. This ordinary differential equation is to be understood in the weak sense when $\psi$ is not twice differentiable. Equation \eqref{EDO} is nothing but a Schrödinger equation with a linear potential.

We shall prove that $\kappa \neq 0$ does not lead to a local minimizer. Indeed, in that case one can check that $y(z)=\psi(\kappa^{-1/3}z-\kappa^{-1}\ell)$ satisfies $y^{\prime\prime}(z)=z\,y(z)$. Therefore, solutions of \eqref{EDO} can be written as $\alpha\, {\rm Ai}(z) + \beta\, {\rm Bi}(z)$ with $z = \kappa^{1/3} x + \kappa^{-2/3}\ell$ and $\alpha, \beta \in \mathbb{R}$, where ${\rm Ai}$ and ${\rm Bi}$ denote the Airy functions of the first and second kinds. It is known (see \cite{DLMF}, for instance) that ${\rm Ai}(z)$ is positive for $z \geq 0$, decreases rapidly to $0$ as $z \to \infty$, whereas ${\rm Bi}(z)$ diverges to $+\infty$. In addition, both functions oscillate around $0$ as $z \to -\infty$, and when one of them vanishes the other alternates between positive and negative values. Hence, in order to ensure that $\psi\geq 0$ and that $\psi^2$ is a probability density, one necessarily has $\kappa>0$ and $\beta=0$. Besides, one has 
\begin{equation}
\alpha=\sqrt{\frac{a}{\int_{b}^{\infty} \operatorname{Ai}^2(x)\,dx}}
\quad\text{where}\quad  
a=\kappa^{1/3}
\quad\text{and}\quad  
b=\kappa^{-2/3}\ell. 
\end{equation}
Then, observe by integration by parts  that
\begin{align}
\int_{0}^\infty (\psi^{\prime}(x))^2 \, dx
&= -\psi(0)\psi^\prime(0) - \int_{0}^\infty (\ell+\kappa x)\,\psi(x)^2 \, dx\label{eq:onycroit} \\
&= \frac{-\operatorname{Ai}(b)\operatorname{Ai}^\prime(b)}{\int_b^\infty \operatorname{Ai}^2(x)\,dx}\, a^3 - \ell - \frac{\kappa}{2\lambda}
= \left(\frac{-\operatorname{Ai}(b)\operatorname{Ai}^\prime(b)}{\int_b^\infty \operatorname{Ai}^2(x)\,dx}-\frac{1}{2\lambda}\right) a^3 - b a^2.\label{eq:jenpeuxplus}
\end{align}
Denote by $c$ the coefficient  in front of $a^3$ in the right-hand side of equation \eqref{eq:jenpeuxplus}. 

If $c=0$, then $b\neq 0$ (otherwise the left-hand side of \eqref{eq:onycroit} would vanish), and thus, keeping $b$ fixed and adjusting $a$ according to the sign of $b$ yields a critical function $\psi$ with lower entropy. In other words, the function $\psi$ corresponding to such $(a,b)$ is not a local minimizer. If $c\neq 0$ and $b=0$, the latter conclusion remains the same for similar reasons. Hence, we can assume that $c\neq 0$ and $b\neq 0$. But the case where $c>0$ and $b>0$ is also not possible if we want a local minimizer. Indeed, in that case the only local minimum of $a\mapsto c a^3- b a^2$ on $(0,\infty)$ is attained at $a=2b/3c$, but at this point the minimum is negative. The case $c>0$ and $b<0$ is also excluded, since then $a\mapsto c a^3- b a^2$ is increasing on $(0,\infty)$. Similarly, we exclude the cases $c<0$ and $b<0$, as well as $c<0$ and $b>0$.

Finally, we necessarily have $\kappa=0$ in order to obtain a local minimizer. It then follows that $\ell=\lambda^2$ and $\psi(x)=\sqrt{2\lambda}\,e^{-\lambda x}$. Using this, we obtain ${\psi^\prime(x)}/{\psi(x)}=-\lambda$. Subsequently, we recover the drift reflected diffusion as given in \eqref{reflect}.

\section{Extended Models}
\label{sec:extended}
\setcounter{equation}{0}

In this section, we discuss how our results can be used to study more general lattices, specifically two-dimensional generalizations and non-nearest-neighbor extensions. 
We note that functional limit theorems for random walks with general increments were recently investigated in~\cite{VoPey}, 
where convergence to skew Brownian motion was established. 
This suggests that similar approaches could be applied in the context of spider networks.

\subsection{The true spider lattice}

\begin{figure}[H]
	\centering
	\includegraphics[scale=0.7]{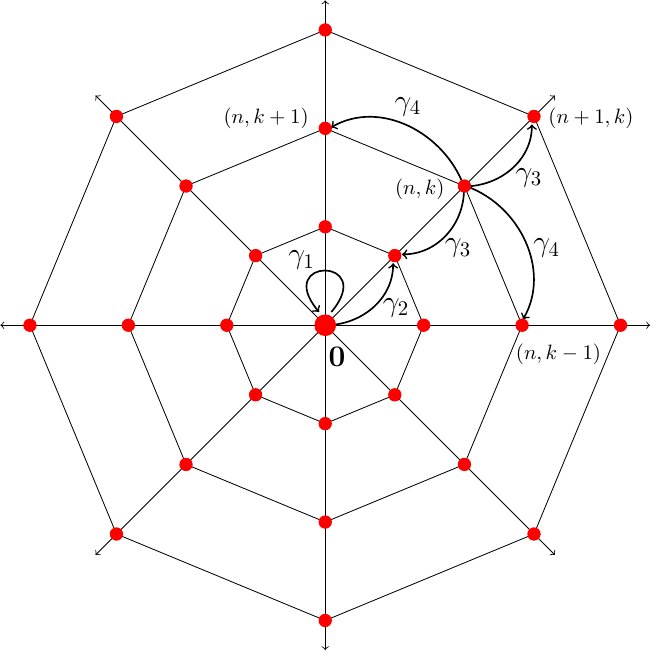}
	\caption{The true spider lattice embbeded into $\mathbb R^2$}
	\label{truespider}
\end{figure}

Consider the spider lattice with $N$ rays, where the $n$th level of each ray is connected to the $n$th level of its two neighboring rays (see Figure \ref{truespider}). In the sequel, we will  assume rotational invariance of the graph. Therefore, the weighted structure can be described by a four-parameter family $\gamma = (\gamma_1, \gamma_2, \gamma_3, \gamma_4)$, with, for all $n \geq 1$ and $k \in \mathbb{Z}/N\mathbb{Z}$, $A(\mathbf{0}, \mathbf{0}) = \gamma_1$, $A(\mathbf{0}; (1, k)) = \gamma_2$, $A((n, k); (n \pm 1, k)) = \gamma_3$, and $A((n, k); (n, k \pm 1)) = \gamma_4$. All other weights being equal to 0.

By symmetry, the reduced model corresponds to the spider lattice in \ref{reduced} with $N = 1$, $\gamma_2 := N \gamma_2$, but with additional loops of weight $2\gamma_4$ above each positive integer. The spectral radius of this simplified model can be computed as before, by replacing the Catalan generating function in (\ref{G_S}) with the Motzkin one, which is given by
\begin{equation}\label{Mot}
S(z) = \frac{1 - 2\gamma_4 z - \sqrt{(1 - 2\gamma_4z)^2 - 4\gamma_3^2 z^2}}{2\gamma_3^2 z^2}.
\end{equation}
Again, equation (\ref{Mot}) can be obtained by applying the ARCH decomposition, leading to
\begin{equation}\label{arch2}
S(z) = \frac{1}{1 - 2\gamma_4 z - \gamma_3^2 z^2 S(z)}.
\end{equation}
Proposition \ref{reduction}, along with the same arguments as in the proof of Proposition \ref{rho}, allows us to conclude that the spectral radius of the true spider lattice satisfies $\rho = 2(\gamma_3 + \gamma_4)$ or
\begin{equation}
\rho=
\dfrac{2\gamma_3(\gamma_1^2+(N\gamma_2)^2)-4N\gamma_1\gamma_2\gamma_4}{\gamma_1(2\gamma_3-N\gamma_2)-2N\gamma_2\gamma_4+N\gamma_2\sqrt{\gamma_1^2+4\gamma_3(N\gamma_2-\gamma_3)+4\gamma_4(\gamma_4-\gamma_1)}},
\end{equation}
according wheither or not $\Lambda:=2(\gamma_3+\gamma_4)-\gamma_1-N\gamma_2\geq 0$.  Besides, Propositions \ref{reg_d}, \ref{att_d}, and \ref{rep_d} can be easily generalized, at least in the $R$-recurrent situation $\Lambda\leq 0$. When $\Lambda>0$, the unique symmetric positive $\rho$-harmonic function is given, for all $n \geq 0$ and $k \in \mathbb{Z}/N\mathbb{Z}$, by
\begin{equation}
\psi(n,k) = 1 + \frac{\Lambda}{N\gamma_2} n,
\end{equation}
However, we point out that it is not certain that all positive harmonic functions are symmetric in that case (see Remark \ref{lose}), further work is needed to determine whether non-symmetric harmonic functions exist. Note that, for a symmetric MERW, the probability of moving from a ray to one of its two neighboring rays is equal to $\gamma_4/\rho$.  Regarding the scaling limits, it seems necessary to let $N \to \infty$ to obtain interesting limits. However, this case is outside the scope of this article and is left for future work.

\subsection{Non-nearest neighbor situation}

\begin{figure}[H]
	\centering
	\includegraphics[scale=0.99]{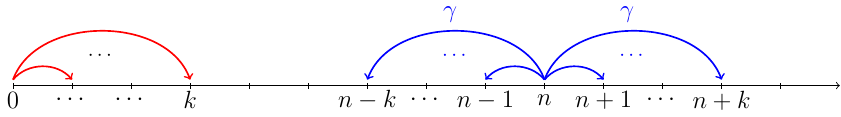}
\caption{A model with $2k$-neighbors}
	\label{nonnearest}
\end{figure}

It is natural to ask whether our results can be extended to the non-nearest neighbor case. We will explore the challenges involved but, for simplicity, we assume $N=1$ (we refer to Figure \ref{nonnearest}). 

Fix $k \geq 1$ and $\gamma > 0$, and assume that $A(n, n \pm i) = \gamma$ for all $n \geq k$ and $1 \leq i \leq k$, and that $A(n, m) = 0$ whenever $|m - n| > k$. Additionally, suppose there are a finite number of non-zero values $A(i, j)$ for $0 \leq i < k$ and $j \geq 0$, and that the corresponding weighted graph is connected.

Here, the arch decomposition used throughout this paper, particularly in (\ref{G}), does not apply as easily, making it unlikely to obtain a necessary and sufficient criterion for $R$-transience. However, we can provide a sufficient criterion. Let us introduce
\begin{equation}
\delta := \max\left\{\sum_{j \geq 0} A(i, j) : 0 \leq i < k\right\}.
\end{equation}

\begin{thm}
	If $\delta \leq 2\gamma k$, then the spectral radius of $A$ is equal to $\rho = 2\gamma k$. If strict inequality holds, then $A$ is $R$-transient. In that case, there exists a unique MERW. More precisely, the unique positive $\rho$-harmonic function can be written as  
	\begin{equation}\label{harmok}
	\psi(n) = a + bn + \mathcal{O}(\beta^n),
	\end{equation}
	for some $a \in \mathbb{R}$, $b > 0$, and $0 < \beta < 1$. Moreover, when $X_0 / L \longrightarrow y \geq 0$, we have
	\begin{equation}\label{scale11}
	\left\{\frac{X_{Lt}}{\sqrt{L}}\right\}_{t \geq 0} \xRightarrow[n \to \infty]{} \left\{Y_{\sigma_k^2 t}^{y}\right\}_{t \geq 0}, \quad \text{with} \quad \sigma_k^2 = \frac{1}{k} \sum_{i=1}^k i^2,
	\end{equation}
	where $Y^y$ is the standard three-dimensional Bessel process defined in (\ref{bess3}).
\end{thm}

\begin{proof} We first show that $A$ is $R$-transient and $\rho = 2\gamma k$. To this end, let $B$ be a weighted adjacency matrix such that $A(i,j) \leq B(i,j)$ for all $i,j \geq 0$, $A(i,j) = B(i,j)$ for all $i \geq k$ and $j \geq 0$, and $\sum_{j \geq 0} B(i,j) = 2\gamma k$ for all $i < k$ (and thus for all $i\geq k$). One can easily note that $\sum_{j \geq 0} B^n(i,j) = (2\gamma k)^n$ and thus $\rho \leq 2\gamma k$. Moreover, $A^{(n)}(k, k)$ is greater than the weighted number of excursions of length $n$ from $k$ to $k$ that remain above $k$. Using \cite[Theorem 3, equation (37), p. 61]{BandFlajo2002},  there exists a constant $C > 0$ such that for $n$ sufficiently large,
	\begin{equation}
	A^{(n)}(k,k) \geq C\frac{(2\gamma k)^n}{n^{3/2}},
	\end{equation}
	showing that $\rho \geq  2\gamma k$. Note that the spectral radius of $B$, as defined above, is also equal to $\rho$. Furthermore, when $\delta < 2\gamma k$, one can choose $B \neq A$, differing only on a finite number of edges, in such a way that \cite[Theorem 5 (b)]{swart2019} implies that $A$ is $R$-transient.

Secondly, we prove the uniqueness of $\psi$ and (\ref{harmok}). Note that any $\rho$-harmonic function satisfies
\begin{equation}\label{harmo}
\forall n \geq k,\quad 2k\psi(n) = \sum_{i=1}^k (\psi(n+i) + \psi(n-i)).
\end{equation}
Introduce $P(X) = X^{2k} + \cdots + X^{k+1} - 2k X^k + X^{k-1} + \cdots + 1$, which is the characteristic polynomial of the corresponding linear recurrence relation. We claim that $1$ is the only root of $P$ on the unit circle $\mathbb{U}$. Indeed, let $\omega \in \mathbb{U}$ be such that $P(\omega) = 0$. We can write
\begin{equation}
1 = 2k|\omega| = \left|\sum_{i=0,i \neq k}^{2k} \omega^i\right| \leq 2k.
\end{equation}
However, equality in the triangle inequality implies that all the  $\omega^i$, for $i \neq k$ and $0 \leq i \leq 2k$, are proportional. It follows that $\omega = 1$. Moreover, the root $1$ has multiplicity 2 since
\begin{equation}
P(X) = (X-1)^2 (1+ \cdots + S_{k-1} X^{k-2} + S_{k} X^{k-1} + S_{k-1} X^k + \cdots  + X^{2k-2}),
\end{equation}
where $S_n = 1 + 2 + \cdots + n$. In addition,  we obtain that $1$ is the only non-negative root of $P$. Furthermore, since $P$ is self-reciprocal (i.e., palindromic), we have that, for $z \neq 0$, $P(z) = 0$ if and only if $P(1/z) = 0$. Let $\mathcal{I}$ and $\mathcal{O}$ be, respectively, the sets of all complex roots of $P$, counted with their multiplicities, strictly inside and outside the unit disc. Note that $\mathcal{I}$ and $\mathcal{O}$ both have cardinality $k - 1$. Any real solution of (\ref{harmo}) can be written as $\psi(n) = a + bn + \psi_{\it i}(n) + \psi_{\it o}(n)$, where $\psi_{\it i}$ and $\psi_{\it o}$ correspond to the roots in $\mathcal{I}$ and $\mathcal{O}$. In particular, $\psi_{\it i}$ and $\psi_{\it o}$ are each characterized by $k - 1$ real parameters. The positivity of solutions to linear recurrence relations is an old and difficult problem. Recently, it was shown in \cite{positivity0} (see Theorem 2) that any non-zero solution with no positive characteristic root of maximal modulus oscillates around zero. As a consequence, we obtain that $\psi_{\it o} \equiv 0$ for any positive solution $\psi$. As a consequence, any positive solution is characterized by $k+1$ real coefficients. Besides, such a solution satisfies the $k+1$ linear boundary equations given by $\psi(0) = 1$ and, for all $0 \leq n < k$,
\begin{equation} 
2k\gamma \psi(n) = \sum_{m \geq 0} A(n,m)\psi(m) = \Phi_n(\psi(0), \cdots, \psi(k-1)),
\end{equation}
where $\Phi_n$ is some linear functional. We deduce that there exists at most one positive solution. By Theorem \ref{pruitt}, in particular Remark \ref{pruitt2}, there exists at least one positive harmonic function, which proves the existence and uniqueness. Finally,  expression (\ref{harmok}) can be simply obtained by choosing $\beta = \sup_{\zeta \in \mathcal{I}} |\zeta| + \delta < 1$ for some  $\delta > 0$.

It remains to prove that $b > 0$. First, if $b = 0$ and $a = 0$, then Theorem 2 in \cite{positivity0} still applies and ensures that the solution is either zero or oscillating, which is not the case. So if $b = 0$, we must have $a \neq 0$. Then, one can check that for all $1\leq i\leq k$,
\begin{equation}
\frac{\psi(n+i) - \psi(n-i)}{\psi(n)} = \mathcal{O}\left(\beta^n\right).
\end{equation}
In other words, the drift of the corresponding MERW goes to zero exponentially fast. However,  it is well known that the resulting Markov chain is recurrent in that case. For instance, by noting that there exists $C > 0$ such that
\begin{equation}
\mathbb{E}[\sqrt{X_{k+1}} - \sqrt{X_k} | X_k = n] \underset{n \to \infty}{\sim} -\frac{C}{n^{3/2}},
\end{equation}
one can apply a standard Foster-Lyapunov method (see \cite[Proposition 2.2]{Comets} for instance) and prove the recurrence. But this contradicts the fact that $A$ is $R$-transient. Finally, we conclude that necessarily $b>0$.

Finally, let us briefly explain how to deduce the scaling limit (\ref{scale11}). First, note that
\begin{equation}
\frac{\psi(n+i)}{\psi(n)} = 1 + \frac{i}{n} + \mathcal{O}\left(\frac{1}{n^2}\right).
\end{equation}
Hence, for all sufficiently smooth functions $f$, we get 
\begin{equation}
\mathbb{E}\left[f\left(\frac{X_{n+1}}{\sqrt{L}}\right) - f\left(\frac{X_n}{\sqrt{L}}\right)\bigg|\frac{X_n}{\sqrt{L}} = x\right] \underset{L \to \infty}{\sim} \frac{\sigma_k^2}{L} \left(\frac{1}{2} f^{\prime\prime}(x) + \frac{1}{x} f^{\prime}(x)\right).
\end{equation}
Adapting the proof  given in Section \ref{sec:proof}, we can obtain the desired scaling limit.
\end{proof}

	The case $\delta > 2\gamma k$ seems out of reach, even for the computation of the spectral radius $\rho$. The latter can be greater than or equal to $2\gamma k$. Moreover, one can choose $\delta$ as large as desired while keeping $\rho = 2\gamma k$, $R$-transience, and the results above. For instance, taking $k = 2$  with  $A(0,1) = \delta$, $A(1,2) = \varepsilon$ and $A(0,i)=A(1,j)=0$ for all $i\neq 1$ and $j\neq 2$,  we obtain by  the Arch decomposition that the generating function of the walks starting and ending at 0 satisfies
	\begin{equation}
	F_{0,0}(z) = \frac{1}{1 - \delta \varepsilon z^3 F_{2,2}^{\geq 1}(z)},
	\end{equation}
	where, the generating function $F_{2,2}^{\geq 1}$ corresponds to the walks starting and ending at 2, remaining in $\{1,2,\dots\}$. As previously, the subgraph $\{1,2,\cdots\}$ is $R$-transient with spectral radius $2\gamma k$. In particular, $F_{2,2}^{\geq 1}(1/(2\gamma k))$ is finite, and thus, for all $\delta > 0$, it is possible to choose $\varepsilon > 0$ such that $F_{0,0}(1/\rho) < \infty$, implying the $R$-transience and $\rho = 2\gamma k$.

\section{Proofs of the Functional Scaling Limits}
\label{sec:proof}
\setcounter{equation}{0}

We will prove these theorems by following the standard approach: first establishing tightness and then identifying the limit. In what follows, let $\mathbb P^L$ represent the law on $\mathbb C$ of the scaled MERW such that under the probability distribution $\mathbb P^L(d\omega)$, the sequence $\left\{{\omega(t)}\right\}_{t\geq 0}$ is distributed as the left-hand side of (\ref{scale1}), (\ref{scale2}), or (\ref{scale3}), depending on the assumptions. The expectation under this probability distribution will be denoted by $\mathbb E^L$.

\subsection{Tightness}

{\it The submartingale argument.} To prove tightness, we will employ a submartingale argument as found in \cite[chap. 1.4.]{ref7}. While this method is applicable to continuous stochastic processes that take their values in $\mathbb R^d$, it can be readily extended to metric spaces, as indicated in \cite[Theorem 2.1.]{FredWentz}. More specifically, we will make use of the following result.

\begin{thm}\label{tight}
	Let $\{\mathbb P^L : L>0\}$ be a family of probability distributions on ${\bf C}$ satisfying
	\begin{equation}
	\mathbb P^L\left(\left\{\omega\in {\bf C}: \forall n\geq 0,\, \omega \text{ is linear over } \left[\frac{n}{L},\frac{n+1}{L}\right]\right\}\right)=1.
	\end{equation}
	Further, assume that for every $\varepsilon>0$, there exist $A_\varepsilon$ and $L_\varepsilon>0$ such that for any $y\in \mathcal G$, there is a function $f_\varepsilon^y $ on $\mathcal G$ satisfying 
	\begin{equation}
	(i)  \; f_\varepsilon^y(y)=1,\quad
	(ii)  \; f_\varepsilon^y(x)=0 \;\text{ if }\; d(x,y)\geq \varepsilon,\quad
	(iii)  \; 0\leq f_\varepsilon^y(x)\leq 1,
	\end{equation}
	and for all $L\geq L_\varepsilon$, 
	\begin{equation}\label{subdis}
	\left\{f_\varepsilon^y\left(\omega\left({\frac{n}{L}}\right)\right)+A_\varepsilon \frac{n}{L}\right\}_{n\geq 0}\; \text{is a $\mathbb P^L(d\omega)$-submartingale.}
	\end{equation}
	Then, as $L$ goes to infinity, the family of probability measures $\{\mathbb P^L : L>0\}$ is tight.
\end{thm}

\noindent
{\it Ito's formula for regular and attractive cases.}
Let $\{X_n\}_{n\geq0}$ be the MERW satisfying the assumptions of Theorem \ref{reg_c} or \ref{att_c}, and let $\overline{X}_n$ denote its first component. Let $\mathbf f$ be a smooth function on $\mathcal G$ with bounded derivatives on each leg $G_i$, $1\leq i\leq N$. For all $x=(\overline{x},i) \in\mathcal G$ with $\overline{x}\geq 1$, recall that $\mathbf f_i(\overline x)=\mathbf f(\overline x,i)$ and define
\begin{multline}\label{laplace}
\Delta_L \mathbf f(x)=\Delta_L \mathbf f_i(\overline{x})=\mathbf f_i\left(\frac{\overline{x}+1}{\sqrt L}\right)-\mathbf f_i\left(\frac{\overline{x}-1}{\sqrt L}\right)\\
\text{and} \quad \Delta^2_L \mathbf f(x)=\Delta^2_L \mathbf f_i(\overline{x})=\mathbf f_i\left(\frac{\overline{x}+1}{\sqrt L}\right)-2\mathbf f_i\left(\frac{\overline{x}}{\sqrt L}\right)+\mathbf f_i\left(\frac{\overline{x}-1}{\sqrt L}\right).
\end{multline}
It is noteworthy that 
$|\Delta_L \mathbf f(x)|\leq {2\|\mathbf f^\prime\|_\infty}/{\sqrt L}$ 
and $|\Delta^2_L \mathbf f(x)|\leq {2\|\mathbf f^{\prime\prime}\|_\infty}/{L}$.
The drift is defined for all $x=(\overline{x},i)\neq \mathbf 0$ by
\begin{equation}\label{drift}
D(x)=D_i(\overline{x})=\mathbb E[\overline{X}_{n+1}-\overline{X}_n| X_n=x].
\end{equation}
Utilizing the classical discrete-time version of Ito's formula, which follows from \cite[p. 132]{ref27} for instance, we can express for all $0\leq m\leq n$,
\begin{multline}\label{itodiscret22}
\mathbf f\left(\frac{X_{n}}{\sqrt L}\right)=\mathbf f\left(\frac{ X_{m}}{\sqrt L}\right)+\frac{1}{2} \sum_{i=1}^{N}  \sum_{k=m}^{n-1}  \left(\Delta^2_L \mathbf f_i(\overline{X}_k)+D_i(\overline{X}_k)\Delta_L \mathbf f_i(\overline{X}_k)\right)\mathds 1_{\{X_k\in  G_i\setminus\{\mathbf 0\}\}}\\
+\#\{m\leq k\leq n-1: X_k=\mathbf 0\}\sum_{i=1}^{N}P_i(0,1)\left(\mathbf f_i\left(\frac{1}{\sqrt L}\right)-\mathbf f(\mathbf 0)\right)+M_{n}-M_{m},
\end{multline}
where $\{M_n\}_{n\geq 0}$ is a square integrable $(\mathcal F_n)$-martingale. Here, $\# E$ denotes the cardinality of a set $E$ and $\mathcal F_n=\sigma(X_0,\cdots, X_n)$.\\

\noindent
{\it Taylor expansions with respect to $\gamma$}. To go further, introduce
\begin{equation}
\mathcal H=  \left\{(x,y_1,\cdots,y_N,z)\in \mathbb R\times \mathbb R^{N}\times\mathbb R :  x+\sum^N_{k=1}y_k-2z=0\right\}.
\end{equation}
For a given $u\in\mathbb R\times \mathbb R^{N}\times\mathbb R$, let $\pi(u)=(\pi_1(u),\pi_2(u),\pi_3(u))$ be the orthogonal projection of $u$ onto $\mathcal H$, and let $\delta(u)=u-\pi(u)=(\delta_1(u),\delta_2(u),\delta_3(u))$. The $i$-th component of $\pi_2(u)$ or $\delta_2(u)$ is denoted by $\pi_2^{(i)}(u)$ or $\delta_2^{(i)}(u)$ for any $1\leq i\leq N$. Observe that for all $1\leq i\leq N$,
\begin{equation}
\delta_1(\gamma)=\delta_2^{(i)}(\gamma)=-\frac{\Lambda}{N+5} \quad\text{and}\quad \delta_3(\gamma)=\frac{2\Lambda}{N+5}.
\end{equation}
Let us set
\begin{equation}
G(x,y,z)=\frac{z(x^2+y^2)}{x(z-y)+y\sqrt{x^2+2yz-z^2}}.
\end{equation}
We now turn to the attractive case $\Lambda<0$. Note that $\rho=G(\gamma_1,\mathcal S_2,2\gamma_3)$. Standard computation shows that $\nabla G(x,y,z)=(0,0,1)$ when $z=x+y$.  
Furthermore, let $\Gamma(\gamma) \equiv \Gamma$ to highlight the dependence of $\Gamma$ in Theorem \ref{att_c} with respect to the parameters $\gamma\in\mathbb R\times\mathbb R^{N}\times \mathbb R$. We get
\begin{equation}
\nabla\Gamma(\pi(\gamma))= \left(\sum_{i=1}^N\pi_2^{(i)}(\gamma)\right)^{-1}\left(-1,\dots,-1,2\right).
\end{equation}
The first order Taylor expansion of $\Gamma(\gamma)$ at $\pi(\gamma)$ in the attractive case becomes
\begin{equation}\label{asympgamma}
\Gamma=1+\frac{\Lambda}{\mathcal S_2+\frac{N\Lambda}{N+5}}+\mathcal O(\Lambda^2)=1-\frac{\lambda}{\sqrt L}+o\left(\frac{1}{\sqrt L}\right).
\end{equation} 
Besides, still in the attractive case, it holds that  
\begin{equation}
D_i(\overline x)=\frac{\gamma_3(\Gamma^2-1)}{\rho\Gamma}=\frac{\Gamma^2-1}{\Gamma^2+1}.
\end{equation}
Consequently, for $\overline x\geq 1$, one has
\begin{equation}\label{driftasymp}
D_i(\overline x)= -\frac{\lambda}{\sqrt L}+o\left(\frac{1}{\sqrt L}\right).
\end{equation}
In the regular case, we have $D_i(\overline x)=0$. In both scenarios, we obtain uniformly on $G\setminus \{\mathbf 0\}$:
\begin{equation}\label{bound}
\left| \frac{1}{2}\left(\Delta^2_L\mathbf f(x)+D_i(\overline x)\Delta_L \mathbf f(x)\right)\right|=\mathcal O\left( \frac{\|\mathbf f^\prime\|_\infty+\|\mathbf f^{\prime\prime}\|_\infty}{L}\right).
\end{equation}
~

\begin{figure}[!h]
	\begin{center}
		\includegraphics[scale=0.7]{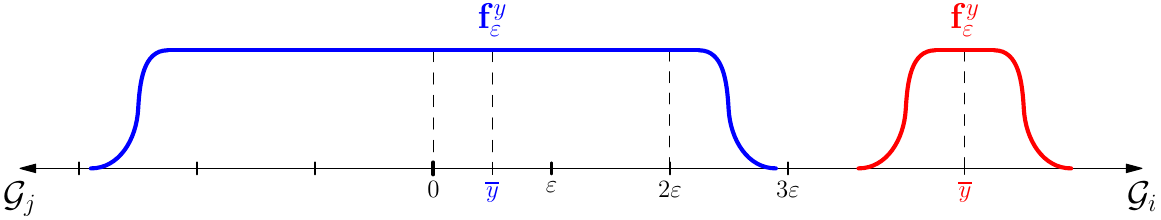}
	\end{center}
	\caption{Lyapunov functions}
	\label{funct}
\end{figure}

\noindent
{\it The Lyapunov functions.} Let $\varepsilon>0$ and let $f_\varepsilon$ be a smooth even function on $\mathbb R$ satisfying the following properties: $0\leq f_\varepsilon\leq 1$, $f_\varepsilon(x)=0$ for $|x|\geq \varepsilon$, $f_\varepsilon(x)=1$ for $|x|\leq \varepsilon/2$, and $f_\varepsilon$ is non-increasing on $[0,\infty)$. Given $y=(\overline{y},i) \in \mathcal G$, we define for all $x=(\overline{x},j) \in \mathcal G$,
\begin{equation}\label{lyapu1}
\mathbf f_\varepsilon^y(x)=
\left\{\begin{array}{ll}
f_\varepsilon(\overline{x}-\overline{y}) &  \text{if $\overline{y}\geq 2\varepsilon$ and $j=i,$}\\[5pt]
f_\varepsilon(\overline{x}-2\varepsilon) &  \text{if $\overline{y}< 2\varepsilon$ and $\overline{x}\geq 2\varepsilon$,}\\[5pt]
1 &  \text{if $\overline{y}< 2\varepsilon$ and $\overline{x}<2\varepsilon$,}\\[5pt]
0 & \text{otherwise.}
\end{array}\right.
\end{equation}
It is evident that $\mathbf f_\varepsilon^y(y)=1$, $0\leq \mathbf f_\varepsilon^y\leq 1$, and $\mathbf f_\varepsilon^y(x)=0$ whenever $d(x,y)\geq 5\varepsilon$. An illustration of this can be found in Figure \ref{funct}. Furthermore, we have $\mathbf f_{\varepsilon,i}^y(1/\sqrt L)-\mathbf f_{\varepsilon,i}^y(\mathbf 0)=0$ for every $L$ satisfying $1/\sqrt L\leq \varepsilon$, given that $\mathbf f_\varepsilon^y$ is flat over the set $\{x\in\mathcal G : d(\mathbf 0,x)\leq \varepsilon\}$. From equations (\ref{bound}) and (\ref{itodiscret22}), we deduce that 
\begin{equation}
\mathbb E\left[\mathbf f_\varepsilon^y\left(\frac{X_{n}}{\sqrt L}\right)-\mathbf f_\varepsilon^y\left(\frac{X_{m}}{\sqrt L}\right)\bigg |\mathcal F_{ m }\right]\geq -A_\varepsilon \left(\frac{n}{L}-\frac{m}{L}\right),
\end{equation}
where $A_\varepsilon$ depends only on the parameters, $\|f_\varepsilon^\prime\|_{\infty}$, $\|f_\varepsilon^{\prime\prime}\|_{\infty}$, and $\varepsilon$. Consequently, we can establish the tightness in Theorems \ref{reg_c} and \ref{att_c} using Theorem \ref{tight}.

\noindent
{\it Focus on the repulsive case.} We only consider the primary distinctions compared to the previous cases. For every $x=(\overline x,i)\neq \mathbf 0$, we have
\begin{equation}\label{driftrep}
D(x)=D_i(\overline x)=\mathbb E\left[\overline X_{n+1}-\overline X_n \big|  X_n= x\right]=\frac{\mu_i\Lambda}{\gamma_2^{(i)}+\mu_i\Lambda\overline x}.
\end{equation}
It can also be verified that
\begin{equation}\label{keypoint}
D_i(\overline x)\Delta_L \mathbf f_{\varepsilon,i}^y(\overline x)=\mathcal O\left(\frac{\|f_\varepsilon^\prime\|}{\varepsilon L}\right)\mathds 1_{\{\overline x\geq \varepsilon \sqrt L\}}.
\end{equation}
Again, we get from the Ito's formula \ref{itodiscret22}:
\begin{equation}
\mathbb E\left[\mathbf f_\varepsilon^y\left(\frac{ X_{n}}{\sqrt L}\right)-\mathbf f_\varepsilon^y\left(\frac{ X_{m}}{\sqrt L}\right)\bigg |\mathcal F_{ m }\right]\geq -A_\varepsilon \left(\frac{n}{L}-\frac{m}{M}\right),
\end{equation}
here for all $L\geq \varepsilon^{-2}$. The constant $A_\varepsilon$ solely depends on the parameters, $\| f_\varepsilon^\prime\|_{\infty}$, $\|f_\varepsilon^{\prime\prime}\|_{\infty}$, and $\varepsilon$. Using again Theorem \ref{tight}, we conclude the proof for tightness.

\subsection{Limit Processes}

Let $\mathbb P^\star$ be a limit point of $\mathbb P^L$ as $L$ approaches infinity, and let $\mathbb E^\star$ represent the associated expectation. We aim to demonstrate that $\mathbb P^\star$ corresponds to the distribution of $\mathbf W^{(\mu,x)}$, $\mathbf Z^{(\mu,x)}$, or $\mathbf Y^{(\mu,x)}$ as specified in Theorems \ref{reg_c}, \ref{att_c}, or \ref{rep_c}, depending on the underlying assumptions. To achieve this, we characterize these stochastic processes in terms of local martingale and/or submartingale problems. Given a continuous function $\mathbf f$ on $\mathcal G$, which is smooth on every ray, we define for all $x=(\overline x,i)\neq \mathbf 0$:
\begin{equation}
\mathcal L_{\mathbf W} \mathbf f(x)=\frac{1}{2}\mathbf f^{\prime\prime}_i(\overline x),\quad \mathcal L_{\mathbf Z} \mathbf f(x)=\frac{1}{2}\mathbf f^{\prime\prime}_i(\overline x)-\lambda \mathbf f^\prime_i(\overline x)
\end{equation}
and 
\begin{equation}
\mathcal L_{\mathbf Y} \mathbf f(x)=\frac{1}{2}\mathbf f^{\prime\prime}_i(\overline x)+\frac{1-\delta_{\mu_i,0}}{\overline x}\mathbf f^\prime_i(\overline x).
\end{equation}
~

\noindent
{\it Martingale and submartingale problems.} The following result is referenced in \cite{Walsh18} and \cite{YanWalsh}.

\begin{thm}\label{martsubpbm}
	Let us define $\mathcal L$ as either $\mathcal L_{\mathbf Z}$ or $\mathcal L_{\mathbf W}$. In the context where $\mathcal L=\mathcal L_{\mathbf Z}$ (resp. $\mathcal L=\mathcal L_{\mathbf W}$), the Walsh diffusion $\mathbf Z^{(\mu,x)}$ (resp. $\mathbf W^{(\mu,x)}$) is  the unique solution $\mathbb Q$ on $\mathbf C$ to the  martingale and submartingale problem described by
	\begin{equation}\label{time}
	\omega_0=x, \quad \int_0^\infty\mathds 1_{\{\omega_s=\mathbf 0\}}ds=0 \quad \text{almost surely with respect to $\mathbb Q(d\omega)$},
	\end{equation}
	and for every sufficiently smooth and bounded function $\mathbf f$ on $\mathcal G$ satisfying either
	\begin{equation}
	i)\; \sum_{i=1}^N \mu_i\mathbf f_i^\prime(0)=0 \quad \text{or} \quad ii)\; \sum_{i=1}^N \mu_i\mathbf f_i^\prime(0)\geq 0,
	\end{equation}
	the stochastic process $\{M^{\mathbf f}_t\}_{t\geq 0}$ is either i) a martingale or ii) a submartingale under the distribution  $\mathbb Q(d\omega)$ where
	\begin{equation}\label{subpbm}
	M^{\mathbf f}_t=\mathbf f(\omega_t)-\int_0^t\mathcal L\mathbf f(\omega_s)\mathds 1_{\{\omega_s\neq \mathbf 0\}}ds 
	\end{equation}
\end{thm}

\begin{rem}
	As a matter of fact, the class of functions $\mathbf f$ considered in \cite{Walsh18,YanWalsh} consists of continuous functions that are twice continuously differentiable on each ray. It is only required that $M^{\mathbf f}$ is a local martingale (resp. local submartingale). By using classical localization and density arguments, one can restrict the domain of functions $\mathbf f$ as stated above and require that $M^{\mathbf f}$ be a martingale (resp. submartingale).   
\end{rem}

Regarding $\mathbf Y^{(\mu,x)}$, the results in \cite{Walsh18,YanWalsh} do not apply. Additionally, the drift of the three-dimensional Bessel process in (\ref{bess3}) is singular at the origin, presenting an additional challenge. However, one can state the following result, the proof of which is provided below.

\begin{thm}\label{newmartpbm}
	The diffusion $\mathbf Y^{(\mu,\mathbf 0)}$ is the unique solution $\mathbb Q$ on $\mathbf C$ of the following martingale problem, given by 
	\begin{equation}\label{time0}
	\omega_0 = \mathbf 0, \quad \int_0^\infty \mathds 1_{\{\omega_s = \mathbf 0\}} ds = 0 \quad \text{$\mathbb Q(d\omega)$-a.s.,}
	\end{equation}
	for all $0 \leq i \leq N$,
	\begin{equation}\label{support}
	\mathbb Q\left(\left\{ \omega \in \mathbf C : \forall t \geq 0, \omega_t \in \mathcal G_i \right\}\right) = \mu_i,
	\end{equation}
	and for all sufficiently smooth and bounded functions $\mathbf f$ with compact support included in $\mathcal G \setminus \{\mathbf 0\}$, 
	\begin{equation}\label{subpbm2}	M^{\mathbf f}=\left\{\mathbf f(\omega_t) - \int_0^t \mathcal L_{\mathbf Y} \mathbf f(\omega_s) ds\right\}_{t \geq 0} \text{is a $\mathbb Q(d\omega)$-martingale.}
	\end{equation}
\end{thm}

\begin{proof} 
	It is clear that the distribution of $\mathbf Y^{(\mu,\mathbf 0)}$ solves the martingale problem. The challenge remains to show that it is the unique solution. Initially, we assume that $N=1$. In this case, $\mathbf Y^{(\mathbf 0,\mu)}=Y^0$ is simply the three-dimensional Bessel process starting from $0$. We denote its distribution by $Q_0$ and its infinitesimal generator by $\mathcal L_Y$. For $\eta>\varepsilon>0$, introduce 
	\begin{equation}
	\sigma_\eta(\omega)=\inf\{t\geq 0 : \omega_t= \eta\}\quad\text{and}\quad \tau_\varepsilon(\omega)=\inf\{t\geq 0 : \omega_{\sigma_\eta+t}=\varepsilon\}.
	\end{equation}
	From our assumptions and using standard localization and approximation arguments, we conclude that
	\begin{equation}\label{subpb3}
	\left\{f(\omega_{\sigma_\eta+t\wedge \tau_\varepsilon})-\int_0^{t\wedge \tau_\varepsilon}\mathcal L_{Y} f(\omega_{\sigma_\eta+s})ds\right\}_{t\geq 0} \text{is a $\mathbb Q(\cdot|\sigma_\eta<\infty)$-martingale},
	\end{equation}
	for all sufficiently smooth functions $f$ on $[0,\infty)$. Letting $\varepsilon\to 0$ and applying standard results on martingale problems and stochastic differential equations, we deduce that $\mathbb Q(\cdot|\sigma_\eta<\infty)$ equals the distribution $Q_\eta$ of a three-dimensional Bessel process starting from $\eta$. Moreover, from (\ref{time0}), we find
	$\mathbb Q(\cdot|\sigma_\eta<\infty)\Longrightarrow \mathbb Q$ as $\eta\downarrow 0$ because 
	\begin{equation}
	\mathbb P\left(\bigcup_{\eta\downarrow 0}\{\sigma_\eta<\infty\}\right)=1.
	\end{equation}
	Since the three-dimensional Bessel process is a Feller Markov process, we have $Q_\eta\Longrightarrow Q_0$ as $\eta\downarrow 0$. Thus, we obtain $\mathbb Q=Q_0$. 	Finally, assuming $N$ is arbitrary, for any $1\leq i\leq N$, introduce the stopping time $\tau_i(\omega)=\inf\{t\geq 0 : \omega_t\notin \mathcal G_i\}$. From (\ref{support}) and (\ref{time0}), we deduce that $\tau_i\sim \mu_i\delta_\infty+(1-\mu_i)\delta_0$. Specifically, $\{\tau_i=\infty\}\in\bigcap_{s>0}\mathcal F_s$. In the sequel, we denote the expectation with respect to the conditional probability $\mathbb Q(\cdot |\tau_i=\infty)$ by $\mathbb E^{(i)}$. For fixed $t\geq s>0$, let $\Upsilon$ be a bounded $\mathcal F_s$-measurable random variable. We obtain
	\begin{equation}\label{subpbm4}
	\mathbb E\left[(M_{t\wedge \tau_i}^{\mathbf f}-M_{s\wedge \tau_i}^{\mathbf f})\mathds 1_{\{\tau_i=\infty\}} \Upsilon\right]=\mu_i\mathbb E^{(i)}\left[(M_{t}^{\mathbf f}-M_{s}^{\mathbf f}) \Upsilon \right]=0.
	\end{equation}
	Note that the latter equality also holds for $s=0$ due to continuity. Hence, $M^{\mathbf f}$ is a $\mathbb Q(\cdot|\tau_i=\infty)$-martingale. Utilizing the result for $N=1$, we conclude that $\mathbb Q$ is the distribution of $\mathbf Y^{(\mu,\mathbf 0)}$.
\end{proof}

\noindent
{\it Bound on the local time at zero.} To ensure that conditions (\ref{time}) or (\ref{time0}) hold for any limit point $\mathbb Q^\star$, we require the following lemma.

\begin{lem}\label{occupation}
	Let $\{X_n\}_{n\geq 0}$ be the MERW specified in Theorems \ref{reg_c}, \ref{att_c}, or \ref{rep_c}. Remember that $\lambda>0$ is provided in Theorem \ref{att_c}. For any $u, v, \eta > 0$, there exists a positive constant $C_{v,\lambda}$ such that 
	\begin{equation}\label{bound2}
	\limsup_{L\to\infty}\mathbb E\left[\frac{\#\{0\leq k\leq \lceil Lu \rceil : d(X_k,\mathbf 0)\leq \eta \sqrt L\}}{L}\right]\leq\frac{(u+v)(1-e^{-{2\lambda}\eta})}{C_{v,\lambda}}. 
	\end{equation}
\end{lem}

\begin{proof} 
	Firstly, using a simple coupling argument, we can reduce the problem to proving this lemma for a MERW starting from the origin under the assumptions of Theorem \ref{att_c}. Indeed, let $Q$ denote the Markov kernel associated with a regular or repulsive MERW. We represent the corresponding parameters by $g\in\mathbb R^{N+2}$. Choose $\zeta\in\mathbb R^{N+2}$ such that $\zeta_1\geq 0$, $\zeta_2^{(k)}\geq 0$ for all $1\leq k\leq N$ with $\sum_{k}\zeta_2^{(k)}>0$, $\zeta_3>0$, and $\zeta_1+\sum_{k}\zeta_2^{(k)}=2\zeta_3$ subject to the condition
	\begin{equation}
	\frac{\zeta_1}{2\zeta_3}> \frac{g_1}{2g_3}.
	\end{equation}
	Let $P$ be the transition kernel associated with the MERW with parameters $\gamma$ as in Theorem \ref{att_c}, converging to $\zeta$. We then have
	\begin{equation}
	P(\mathbf 0,\mathbf 0)=\frac{\zeta_1}{2\zeta_3}+\mathcal O\left(\frac{1}{\sqrt L}\right).
	\end{equation}
	From this, we get $P(\mathbf 0,\mathbf 0)\geq Q(\mathbf 0,\mathbf 0)$ for large $L$. Also, noting that $P_k(n,n-1)\geq Q_k(n,n-1)$ for all $n\geq 1$ and $1\leq k\leq N$, constructing a coupling for which $X^{P}\leq X^{Q}$ becomes feasible, where $X^{P}$ and $X^{Q}$ are the MERW associated with $P$ and $Q$. 
	
	Secondly, let $\pi$ be the invariant probability distribution of the MERW $\{X_n\}_{n\geq 0}$ in Theorem \ref{att_c}. Using (\ref{asympgamma}) and (\ref{inv}), we deduce
	\begin{equation}\label{expo}
	\sum_{k=1}^N\pi_k((\eta\sqrt L,\infty))=\frac{\mathcal S_2{}\Gamma^{2\lfloor \eta\sqrt L\rfloor}}{\gamma_3(1-\Gamma^2)+\mathcal S_2\Gamma^2}=e^{-2\lambda\eta}+o\left({1}\right),
	\end{equation} 
	leading to
	\begin{equation}
	\frac{\mathbb E_\pi\left[{\#\{0\leq k<\lceil L(u+v) \rceil : \overline{X}_k\leq \eta\sqrt L\}}\right]}{L}\underset{L\to\infty}{\sim} (u+v)(1-e^{-2\lambda\eta}).
	\end{equation}
	Define $T=\inf\{k\geq 0 : \overline{X}_k=0\}$. Employing the strong Markov property, we have
	\begin{eqnarray}
	\mathbb E_\pi\left[\sum_{k=0}^{\lceil L(u+v) \rceil -1} \mathds 1_{\{ \overline{X}_k\leq \eta \sqrt L\}}\right] & \geq & 
	\sum_{i=0}^{\lceil L(u+v) \rceil -1}\mathbb P_\pi(T=i)\,\mathbb E_0\left[ \sum_{k=0}^{\lceil L(u+v) \rceil -1-i} \mathds 1_{\{ \overline{X}_k\leq \eta\sqrt L\}}\right]\\
	&\geq & \mathbb P_\pi(T< \lceil Lv \rceil) \,\mathbb E_0\left[ \sum_{k=0}^{\lceil Lu \rceil -1} \mathds 1_{\{ \overline{X}_k\leq \eta \sqrt L\}}\right].
	\end{eqnarray}
	Let $\Xi_k$ be a random variable distributed as $\pi$ conditionally on $\{(\overline x,k) : \overline x\geq 1\}$ and let $(\xi_i)_{i\geq 1}$ be a sequence of i.i.d.\@ Rademacher random variables with parameter $p=\gamma_3\Gamma\rho^{-1}$ independent of $\Xi_k$. Set $S_n=\xi_1+\cdots+\xi_n$. Through a simple coupling argument, we can infer
	\begin{equation}
	\mathbb P_\pi(T< \lceil Lv \rceil)\geq \pi(\mathbf 0)+\sum_{k=1}^N \pi_k((0,\infty))\mathbb P\left(S_{\lceil Lv \rceil -1}\leq -\Xi_k\right).
	\end{equation}
	Additionally, we have
	\begin{equation}
	\mathbb P\left(S_{\lceil Lv \rceil -1}\leq -\Xi_k\right)=\mathbb P\left(\frac{S_{\lceil Lv \rceil -1}+\lambda\sqrt{L}v}{\sqrt{Lv}} \leq \frac{-\Xi_k}{\sqrt{Lv}}+\lambda \sqrt v\right).
	\end{equation}
	From (\ref{driftasymp}), we observe that $\mathbb E[S_{\lceil Lv \rceil -1}] = -\lambda\sqrt L v + o(\sqrt L)$ as $L$ approaches infinity. Then, analogously to (\ref{expo}) and with the aid of the central limit theorem, we deduce
	\begin{equation}
	\frac{S_{\lceil Lv \rceil -1}+\lambda\sqrt{L}v}{\sqrt{Lv}}\otimes \frac{\Xi_k}{\sqrt {L}} \xRightarrow[L\to\infty]{} (U,V)\sim \mathcal N(0,1)\otimes \mathcal E(2\lambda).
	\end{equation}
	Given that $\pi(\mathbf 0)$ converges to $0$, we conclude
	\begin{equation}
	\liminf_{L\to\infty}\mathbb P_\pi(T< \lceil Lv \rceil) \geq \mathbb P(U+\sqrt{2} V\leq \lambda\sqrt{v}).
	\end{equation}
	Lastly, combining this with the aforementioned equations, the proof is completed.
\end{proof}

\noindent
{\it Identification of the limit.} To go further, one can express 
\begin{equation}\label{itodiscretoop}
\mathbf f\left(\omega\left(\frac{n}{L}\right)\right)=\mathbf f\left(\omega\left(\frac{m}{L}\right)\right)
+\frac{1}{L}\sum_{k=m}^{n-1}  \mathcal L_L\mathbf f\left(\omega\left(\frac{k}{L}\right)\right)+M_{n}-M_m,
\end{equation}
where $M$ is a square-integrable $\mathbb P^L$-martingale. For every $x=(\overline x,i)\in \frac{1}{\sqrt L}\cdot G\setminus\{\mathbf 0\}$,
\begin{multline}\label{operator}
\mathcal L_L\mathbf f(x)=\frac{L}{2}\left(\mathbf f_i\left(\overline x+\frac{1}{\sqrt L}\right)+\mathbf f_i\left(\overline x-\frac{1}{\sqrt L}\right)-2\mathbf f_i(\overline x)\right)\\
+\frac{L D_i\left(\overline x\sqrt{L}\right)}{2}\left(\mathbf f_i\left(\overline x+\frac{1}{\sqrt L}\right)-\mathbf f_i\left(\overline x-\frac{1}{\sqrt L}\right)\right),
\end{multline}
and 
\begin{equation}\label{A20}
\mathcal L_L\mathbf f(\mathbf 0)=L\sum_{i=1}^N P_i(0,1)\left(\mathbf f_i\left(\frac{1}{\sqrt L}\right)-\mathbf f(\mathbf 0)\right).
\end{equation}
Subsequently, we assume that the test functions $\mathbf f$ are sufficiently smooth and bounded, along with their derivatives, on each ray. Importantly, uniformly on $\mathbb R$, we find that
\begin{equation}\label{A2}
\frac{L}{2}\left(\mathbf f_i\left(\overline x+\frac{1}{\sqrt L}\right)+\mathbf f_i\left(\overline x-\frac{1}{\sqrt L}\right)-2\mathbf f_i(\overline x)\right)=\frac{1}{2}\mathbf f^{\prime\prime}_i(\overline x)+\mathcal O\left(\frac{1}{\sqrt L}\right).
\end{equation} 
It is noteworthy that the constant in the big $\mathcal O$ depends exclusively on $\|\mathbf f_i^{\prime\prime\prime}\|_{\infty}$. Our analysis will center on the remaining terms in (\ref{operator}) and (\ref{A20}).\\

\noindent
{\it {Focus on the} regular and attractive cases.} We shall prove that $\mathbb P^\ast$ is {the} solution of the well-posed martingale/submartingale problem in {Theorem} \ref{martsubpbm}. First, the assumptions of Theorems \ref{reg_c} and \ref{att_c} allow us to see that for some $\alpha>0$, one has
\begin{equation}\label{expL0}
\mathcal L_L\mathbf f(\mathbf 0)=\alpha
\sum_{i=1}^N \mu_i\mathbf f^\prime_i(0)\sqrt L+\mathcal O(1).
\end{equation}
Here we use $P_i( 0,1)=\frac{\zeta_2^{(i)}}{2\zeta_3}+ O\left(\frac{1}{\sqrt{L}}\right)$ for the attractive case. Recall that $D_i(\overline x\sqrt L)=0$ as soon as  $\overline x \in \frac{1}{\sqrt L}\cdot \mathbb N$ in the regular case whereas in the attractive case we obtain from (\ref{driftasymp}) that 
\begin{equation}\label{keypoint2}
\frac{L D_i\left(\overline x\sqrt{L}\right)}{2}\left(\mathbf f_i\left(\overline x+\frac{1}{\sqrt L}\right)-\mathbf f_i\left(\overline x-\frac{1}{\sqrt L}\right)\right)=-\lambda\,\mathbf f_i^\prime(\overline x)+o(1),
\end{equation}
By using (\ref{A2}) we deduce for $\mathcal L\in \{\mathcal L_{\mathbf W},\mathcal L_{\mathbf Z}\}$ according {to} the assumptions that
\begin{equation}\label{opregu}
\mathcal L_L \mathbf f(x)=\left(\mathcal L_{}\mathbf f(x)+o\left(1\right)\right)\mathds 1_{\{x\neq \mathbf 0\}}+(A_L\mathbf f+\mathcal O(1))\mathds 1_{\{x=\mathbf 0\}},
\end{equation}
where $A_L \mathbf f$ satisfies   
\begin{equation}\label{subcond}
\sum_{i=1}^N \mu_i\mathbf f^\prime_i(0)= 0\; \left(\text{resp.}\; \sum_{i=1}^N \mu_i\mathbf f^\prime_i(0)\geq 0\right)\;\Longrightarrow \; A_L\mathbf f= 0 \;\left(\text{resp.}\; A_L\mathbf f\geq 0\right).
\end{equation}
Let $T,\varepsilon,R\geq 0$ be {given} and set for all $\omega\in {\mathbf C}$, 
\begin{equation}
\delta_\omega(T,\varepsilon)=1\wedge \sup \{d(\omega(t),\omega(s)) : |t-s|\leq \varepsilon, 0\leq s,t\leq T\}.
\end{equation}
Here we denote $a\wedge b=\min(a,b)$.
By using (\ref{itodiscretoop}) one can write for all $0\leq s\leq t\leq T$,
\begin{equation}\label{itodiscret4}
\mathbf f\left(\omega(t)\right)=\mathbf f\left(\omega(s)\right)+\frac{1}{L}\sum_{k={\lfloor Ls\rfloor}}^{{\lfloor Lt\rfloor}-1}  
\mathcal L_L\mathbf f\left(\omega\left(\frac{k}{L}\right)\right)+M_{{\lfloor Lt\rfloor}}-M_{{\lfloor Ls\rfloor}}+
\mathcal O\left(\delta_{\omega}\left(T,\frac{1}{L}\right)\right),
\end{equation}
where the constant in the big $\mathcal O$ {depends} only on the $\|\mathbf f^\prime_i\|_{\infty}$ for $1\leq i\leq N$. Furthermore, we get from (\ref{opregu}) that for all $0\leq s\leq t\leq T$,
\begin{multline}
\frac{1}{L}\sum_{k={\lfloor Ls\rfloor}}^{{\lfloor L t\rfloor}-1}  \mathcal L_L\mathbf f\left(\omega\left(\frac{k}{L}\right)\right)\mathds 1_{\left\{\omega\left(\frac{k}{L}\right)\neq\mathbf  0\right\}}  =  \int_{\frac{{\lfloor Ls\rfloor}}{L}}^{\frac{{\lfloor Lt\rfloor}}{L}}\mathcal L_L\mathbf f\left(\omega\left(\frac{{\lfloor Lu\rfloor}}{L}\right)\right)\mathds 1_{\left\{\omega\left(\frac{{\lfloor Lu\rfloor}}{L}\right)\neq \mathbf 0\right\}} du\\
=   \int_s^{t}\mathcal L\mathbf f\left(\omega\left(u\right)\right)\mathds 1_{\left\{\omega\left(u\right)\neq \mathbf 0\right\}} du + o(1)+\mathcal O\left(\delta_\omega\left(T,\frac{1}{L}\right)+\frac{1}{L}\right).
\end{multline}

At this point we need to note that if we assume that ${\lfloor Lt\rfloor}\geq {\lfloor Ls\rfloor}+1$ and $Ls\neq {\lfloor Ls\rfloor}$ then by using  (\ref{expL0}) and (\ref{itodiscretoop}) one has
\begin{equation}
\mathbb E^L\left[M_{{\lfloor Lt\rfloor}}-M_{{\lfloor Ls\rfloor}}|\mathcal F_s\right](\omega)=M_{{\lfloor Ls\rfloor}+1}(\omega)-M_{{\lfloor Ls\rfloor}}(\omega)=\mathcal O\left({\delta_\omega}\left(T,\frac{1}{L}\right)+\frac{1}{\sqrt L}\right).
\end{equation}
Otherwise $\mathbb E^L\left[M_{{\lfloor Lt\rfloor}}-M_{{\lfloor Ls\rfloor}}|\mathcal F_s\right]=0$.
Finally, we deduce that for all $\mathbf f$ sufficiently smooth with bounded derivatives on each ray satisfying the left-hand-side of (\ref{subcond}) one has 
\begin{multline}\label{submartfin}
\mathbb E^L\left[\mathbf f(\omega(t))-\mathbf f(\omega(s))-\int_s^{t}\mathcal L\mathbf f\left(\omega\left(u\right)\right)\mathds 1_{\left\{\omega\left(u\right)\neq \mathbf 0\right\}} du\bigg|\mathcal F_s\right]\geq (\text{resp.}\; =)\\
o(1)+\mathcal O\left(\mathbb E^L\left[{\delta_\omega}\left(T,\frac{1}{L}\right)\right]\right)+\mathcal O\left(\frac{\mathbb E^L[\#\{0\leq k\leq {\lfloor LT\rfloor}: X_k=\mathbf 0\}]}{L}\right).
\end{multline}
Note that the functional into the expectation of the left-hand-side of (\ref{submartfin}) is continuous and bounded with respect to $\omega\in{\mathbf C}$. Besides, it follows from the tightness and Lemma \ref{occupation} (by letting $\eta\to 0$) that the second and the third term in the right-hand-side of  (\ref{submartfin}) goes to $0$ as $L$ goes to infinity. We deduce (\ref{subpbm}). To conclude, it remains to prove (\ref{time}). To this end, one can write  for all $\eta>0$,
\begin{equation}\label{occuplim}
\mathbb E^\star\left[\int_0^T\mathds 1_{\{d(\omega(s),\mathbf 0)< \eta\}}ds\right]
\leq  \liminf_{L\to\infty} \frac{\mathbb E^L[\#\{0\leq k\leq {\lfloor LT\rfloor}: d(X_k,\mathbf 0)<\eta\sqrt L\}]}{L}.
\end{equation}
Here we use $\{\omega\in {\mathbf C} : d(\omega(s),{\mathbf 0})<\eta\}$ is open and the Fatou's Lemma. Applying again Lemma \ref{occupation} and letting $\eta\to 0$ we obtain (\ref{time}). This completes the proof of Theorems \ref{reg_c} and \ref{att_c}.\\

\noindent
{\it Focus on the repulsive case.} The proof follows the main lines as in the regular and attractive cases above and most of the previous notations are kept. For instance, as for the regular and attractive cases, we deduce from Lemma \ref{occupation} that 
\begin{equation}
\mathbb E^\star\left[\int_0^\infty\mathds 1_{\{\omega(s)=\mathbf 0\}}ds\right]=0.
\end{equation}
The test functions $\mathbf f$ we consider are supposed to have a compact support in $\mathcal G\setminus\{\mathbf 0\}$ as in Theorem \ref{newmartpbm}. In particular, assuming $\mu_i\neq 0$,  the asymptotic (\ref{keypoint2}) becomes 
\begin{equation}\label{keypoint3}
\frac{L D_i\left(\overline x\sqrt{L}\right)}{2}\left(\mathbf f_i\left(\overline x+\frac{1}{\sqrt L}\right)-\mathbf f_i\left(\overline x-\frac{1}{\sqrt L}\right)\right)=\frac{\mathbf f_i^\prime(\overline x)}{\overline x}+o(1).
\end{equation}
Here we use (\ref{driftrep}). Furthermore, we need to distinguish whether or not $x=\mathbf 0$ and when $x=(\overline x,i)\neq \mathbf 0$ whether or not $\mu_i=0$.

{{\textbf{1)}} Assume that $x=(\overline x,i)\neq \mathbf 0$ and $\mu_i\neq 0$. One can prove as previously  that 
	\begin{equation}\label{subpbm2bis}
	\left\{\mathbf f(\omega_{t})-\int_0^{t}\mathcal L_{\mathbf Y}\mathbf f(\omega_s)ds\right\}_{t\geq 0} \text{is a } \text{$\mathbb P^\star(d\omega)$-martingale.}
	\end{equation}
	As in the proof of Theorem \ref{newmartpbm}, we deduce that $\mathbb P^\ast\sim \mathbf Y^{(\mu,x)}$ since the hitting time of $0$ of a three-dimensional Bessel process starting from $\overline x>0$ is infinite almost-surely.

	{{\textbf{2)}} Assume that $x=(\overline x,i)\neq \mathbf 0$ and $\mu_i=0$.	Again one has (\ref{subpbm2bis}) and we deduce that the restriction of $\mathbb P^\star$ to the $\sigma$-algebra $\mathcal F_{\tau_0}$ is a standard Brownian motion on the $i$th ray starting to $x$ up to the hitting time of $\mathbf 0$. Then by using the Markov property, we deduce that $\mathbb P^\ast\sim \mathbf Y^{(\mu,x)}$ provided the result is proved assuming $x=\mathbf 0$.

		{{\textbf{3)}} Assume that $x=\mathbf 0$. Once again (\ref{subpbm2bis}) still holds and in order to apply Theorem \ref{newmartpbm} we only need to show that $\mathbb P^\ast$ satisfies (\ref{support}). As a matter of facts,  it is a simple consequence of (\ref{martinb}) and thus $\mathbb P^\ast\sim \mathbf Y^{(\mu,x)}$.
			
			This ends the proof the scaling limits.\hfill$\Box$\\

\noindent
{\bf Acknowledgements\;} The authors are grateful to the Referees and the Associated Editor for their careful reading and valuable comments and remarks, which have significantly improved the manuscript.

\bibliography{biblio.bib}
\bibliographystyle{unsrt}

\end{document}